\documentclass[12pt, a4paper]{report}

\usepackage{amsmath, amssymb, amsthm}
\usepackage[compat2]{geometry}
\usepackage[matrix,arrow]{xy}

\usepackage{xspace}

\newcommand{\eps}{\varepsilon}
\renewcommand{\labelenumi}{\theenumi)}

\newcommand{\comment}[1]{}
\comment{
This is an English translation of the author's Ph.D. thesis, accumulating
his results on a construction of Cohen-Macaulay modules over a polynomial
ring that appeared in the study of Cauchy-Fueter equations. This
construction is generalized from quaternions to arbitrary
finite-dimensional associative algebras. We show that for maximally central
algebras (as introduced by Azumaya) this construction produces
Cohen-Macaulay modules and is an exact functor (tensoring with a bimodule,
actually) and this class of algebras cannot be enlarged. For this class
several invariants of the resulting modules are calculated via a fairly
explicit description of their graded minimal free resolution, that is
constructed from the Eagon-Northcott complex. These results have been
published in Russ. Math. Surveys and Sbornik: Mathematics but for some
proofs, a concise and complete exposition is presented here.
}

\author{O. N. Popov\thanks{141980, Dubna, Moscow Region, Moskovskaya str., 
12, apt. 236. E-mail: \tt popov@mccme.ru.}}
\title{Modules over a Polynomial Ring Obtained from Representations
of a Finite-dimensional Associative Algebra} 
\date{MSC 2000: 13C14, 16P10, 16E30. Keywords: Cauchy--Fueter equations,
Cohen-Macaulay modules, Eagon-Northcott complex, maximally central
(Azumaya) algebras, tensor product.}

\addtocounter{chapter}{-1}

\newcommand{\set}[1]{\mathbb{#1}}
\newcommand{\K}{\Bbbk}
\newcommand{\M}{{\cal M}}
\newcommand{\KK}{{\cal K}}
\newcommand{\mm}{\mathfrak{m}}
\newcommand{\pd}{\mathop{\mathrm{pd}}}
\newcommand{\depth}{\mathop{\mathrm{depth}}\nolimits}
\newcommand{\hht}{\mathop{\mathrm{ht}}}%\ht --- это какой-то примитив
\newcommand{\Id}{{\rm Id}}
\newcommand{\T}{\widetilde}
\newcommand{\Fl}[1]{F_l(#1)}% функтор номер l
\newcommand{\Ff}[1]{F_1(#1)}% функтор номер 1
\newcommand{\Fll}[1]{F_{l'}(#1)}% функтор номер l'

\newcommand{\B}{\overline}
\newcommand{\co}{\colon}
\renewcommand{\geq}{\geqslant}
\renewcommand{\leq}{\leqslant}
\newcommand{\rad}{\mathop{\mathrm{rad}}}
\newcommand{\hence}{\ensuremath{\Rightarrow}}
\newcommand{\xto}{\xrightarrow}

\newcommand{\End}{\mathop{\mathrm{End}}\nolimits}
\newcommand{\Hom}{\mathop{\mathrm{Hom}}\nolimits}
\newcommand{\Tor}{\mathop{\mathrm{Tor}}\nolimits}
\newcommand{\Br}{\mathop{\mathrm{Br}}}
\newcommand{\Ext}{\mathop{\mathrm{Ext}}\nolimits}

\newcommand{\BF}{\mathbb F}
\newcommand{\BZ}{\mathbb Z}
\newcommand{\fp}{\mathfrak p}

\newcommand{\Proof}{\noindent{\bf Proof. }}
\newcommand{\Pf}{\noindent{\bf Proof. }}
\newcommand{\Homgr}{\mathop{\mathrm{Homgr}}^0_R}
\newcommand{\ann}{\mathop{\mathrm{ann}}\nolimits}
\newcommand{\Ker}{\mathop{\mathrm{Ker}}}
\newcommand{\Coker}{\mathop{\mathrm{Coker}}}
\newcommand{\Ass}{\mathop{\mathrm{Ass}}}

\DeclareMathOperator{\rk}{rk}

\newcommand{\myto}[1]{\xrightarrow{#1}}
\newcommand{\Idarrow}{\myto{\Id\cdot}}

\newcommand{\CC}{{\cal C}}

\newtheorem{stat}{Theorem}

\newtheorem{lemma}{Lemma}[chapter]
\newtheorem{satz}[lemma]{Proposition}
\newcommand{\rem}%
{\addvspace{.5\baselineskip}\noindent\refstepcounter{lemma}%
{\bf Remark \thelemma. }} 

\newcounter{definition}[section]
\renewcommand{\thedefinition}{\thesection.\arabic{definition}}
\newcommand{\mydef}%
{\addvspace{.5\baselineskip}\noindent\refstepcounter{definition}%
{\bf Definition \thedefinition.} }

\newtheorem{prop}[definition]{Proposition}
\newtheorem{lm}[definition]{Lemma}
\newcommand{\remm}%
{\addvspace{.5\baselineskip}\noindent\refstepcounter{definition}%
{\bf Remark \thedefinition. }}

\renewcommand{\geq}{\geqslant}
\renewcommand{\leq}{\leqslant}

\newcommand{\wrt}{w.\:r.\:t.\@\xspace}
\newcommand{\eg}{e.\:g.\@\xspace}
\newcommand{\ie}{i.\:e.\@\xspace}
\newcommand{\degr}{$^\circ$}

\begin{document}
\maketitle

\addcontentsline{toc}{chapter}{Contents}
\tableofcontents

\chapter*{Introduction}
\addcontentsline{toc}{chapter}{Introduction}

This paper, concerned with both commutative and noncommutative algebra,
is devoted to studying an algebraic generalization of a problem which
arose from applications of commutative algebra to linear PDEs with constant
coefficients. These applications are rather old: an invariant way to 
describe a system of linear PDEs
$$
\begin{cases}
D_{11} f_1+\dots+D_{1m}f_m=0,&\\
\dotfill &\\
D_{n1} f_1+\dots+D_{nm} f_m=0,&\\
\end{cases}
$$ 
 $f_j$ being functions of several real variables,
$D_{ij}$ differential operators, is to consider a (left) module 
$M$ over the ring of differential operators (a $D$-module) which is the 
quotient of a free module of rank $m$ modulo the submodule generated by 
the rows of the matrix $(D_{ij}).$ Then if we consider the ring of smooth 
functions (analytic functions, distributions) $O$ as a module over the 
ring of differential operators, we easily see that the space of smooth 
(resp., analytical, generalized) solutions of our system can be identified 
with the space of homomorphisms of $D$-modules $\Hom(M,O)$: the generators 
of $M$ are taken to functions $f_j,$ which satisfy the equations, because 
there are relations between the generators. But the ring of differential 
operators with constant coefficients is the ring of (commutative) 
polynomials in operators $\frac{\partial}{\partial x_i}$, partial 
derivations \wrt the variables, so in case of constant coefficients we 
obtain a module over a polynomial ring. Many meaningful properties of the
solutions of a system of differential equations have a natural restatement
in terms of commutative algebra of this module, see, \eg, \cite{Pal}.

The module corresponding this way to Cauchy--Fueter equations defining
quaternionic-differenitiable functions has appeared in~\cite{palam}. This 
is $\M_n=R^4/\langle A_n\rangle,$ where $A_n$ is the matrix $U_1|\dots|U_n,
\ \langle A_n\rangle$ the submodule generated by its columns,
$$U_i=\left(
\begin{array}{rrrr}
 x_i& y_i& z_i& t_i\\
 -y_i& x_i &t_i &-z_i \\
 -z_i&-t_i & x_i &y_i\\
 -t_i&z_i&-y_i & x_i
\end{array}
\right),$$
and $ R=\K[\{x_i,y_i,z_i,t_i\}_{i=1}^n]$ for a field $\K.$ 

The authors showed that the projective dimension of the module equals
$2n-1,$ whence they derived that the flabby dimension of the sheaf
$\mathcal{O}$ of quaternionic-differentiable functions in $n$ variables
also equals $2n-1$ \cite[Theorem~3.1]{palam} and that for any open set
$U\subset\set{H}^n$ $H^i(U,\mathcal{O})=0$ for $i\geq 2n-1$
\cite[Cor.~3.4]{palam}.
%Эти 2 факта обобщали уже известные
%свойства пучка голоморфных функций, в которых на месте $2n-1$ фигурировало
%число $n.$

The authors used some concepts and methods of commutative algebra that we 
are going to recall. See Section~\ref{dpcm} for further details.

{\bf Definition.} Let $R$ be a commutative Noetherian ring, $M$ an
$R$-module. The sequence $a_1,\dots,a_n\in R$ is called \emph{$M$-regular,}
if $(a_1,\dots,a_n)M\ne M$ and for $i$ between $1$ and $n$ the
multiplication by $a_i$ is injective on $M/(a_1,\dots,a_{i-1})M$.

{\bf Definition.} Let $R$ be a commutative Noetherian ring, $M$ an
$R$-module, $I\subset R$ an ideal and $IM\ne M.$ The length $\depth(I,M)$ 
of any maximal $M$-regular sequence in $I$ is called \emph{depth} of $I$ 
on $M.$ When considering the depth of the ideal of the polynomials
vanishing at the origin on graded modules over a polynomial ring, we shall
omit the ideal and talk of the depth of a module.

{\bf Auslander--Buchsbaum formula.} For a graded module $M$ over a 
polynomial ring $R$ one has $\depth M+\pd M=\dim R.$

To compute the projective dimension, the authors used the
Auslander--Buchsbaum formula, not the explicit
construction of the resolution, as they thought the latter too complicated.  
The depth of the module, which is needed to apply the 
Auslander--Buchsbaum formula and equals $2n+1,$ was caclulated in~\cite{palam}
by means of an explicit construction of an $\M_n$-sequence, 
using Gr\"obner bases. The Krull dimension of $\M_n$ was also determined in
the paper, to which end the tangential space to the support of the module
in $\set{C}^{4n}=\mathop{\rm Specm}R$ for $\K=\set{C}$ was considered.
This way the Cohen-Macaulayness of the module, \ie the equality of the Krull
dimension to the depth (Def.~\ref{CM}), was proved.
 
In a subsequent paper~\cite{AL} the authors continued to investigate this
module with Gr\"obner bases, finding the (graded) Betti numbers (\ie the
ranks and degrees of generators for the components of the graded minimal
free resolution of $\M_n$), the Hilbert series (\ie the dimensions of the
homogeneous components of the module) and the multiplicity of $\M_n$ (\ie
the asymptotics of the dimension of a homogeneous component of the module
as a function in its degree). We recall precise definitions and basic
properties of these concepts in Section~\ref{inv}.

Similar studies of other systems of differential operators with constant 
coefficients were undertaken in \cite{Maxw}, \cite{Moisil}, \cite{Dirac}.

As explained above, in~\cite{palam} the matrix $A_n$ was obtained from the
system of linear partial differential equations for
quaternionic-differentiable functions by transposing and replacing the
partial derivatives \wrt different variables by the variables themselves.
Now one can see (cf.~\cite[Introduction]{AL}) that $U_i$ is the matrix of
the left multiplication by  $x_i-y_ii-z_ij-t_ik$ \wrt the basis $1, i, j,
k.$ This interpretation of $A_n$ allows, as E.~S.~Golod has observed, to
understand completely the structure of the module $\M_n,$ in particular
that of its projective resolution. Let us complexify the algebra of
quaternions. As a change of basis in the algebra leads to a conjugation of
the matrix $U_i$ together with a linear change of variables in it, this
change does not affect the structure of $\M_n.$ Therefore the isomorphism
${\set{H}\otimes_{\set{R}} \set{C}}\cong M_2(\set{C})$ ($M_2(\set{C})$
being the matrix algebra) allows, by choosing the basis consisting of
matrix elements, to transform $U_i$ into
$$\newcommand{\Gen}{
\begin{array}{ll}
a_i&b_i\\
c_i&d_i
\end{array}}
\left(
\begin{array}{cc}
\Gen & \mbox{\Large 0}\\
\mbox{\Large 0} & \Gen
\end{array}
\right),
$$
\ie $\M_n=\M'_n\oplus \M'_n,$  $\M'_n$ being the quotient module of $R^2$
modulo the columns of a generic $2\times 2n$-matrix. Now $\M'_n$ is
known~\cite{BV} to be a Cohen-Macaulay module of projective dimension
$2n-2+1=2n-1$ and to have an Eagon-Northcott complex (the Buchsbaum-Rim
complex of~\cite{Ei}, see Section~\ref{WMCE}) for its minimal resolution.
This description of the resolution allows to simplify the proofs of the
main results of~\cite{AL}, see Chapter 3 below.

From these considerations the following generalization of the problem has
emerged. Let $A$ be a finite-dimensional associative algebra with identity over
a field $\K$ with a $\K$-basis $f_1,\dots,f_d$ and $\rho\co A\to M_n(\K)$ be
its matrix representation corresponding to an $A$-module $M,$ $\dim_{\K} M=n.$
Let us fix a positive integer $l$ and consider the polynomial ring
$R=\K[x_{11},\dots,x_{dl}]$ and the module $\Fl M$ over it ($F^R_l(M)$, if
the ring should be mentioned explicitly, as we shall sometimes consider a
polynomial ring with some additional variables), which is the quotient of a
free $R$-module $R^n$ modulo the submodule generated by the columns of the
matrices $\Id_j=\sum_i \rho(f_i)x_{ij},$ $j=1,\dots,l$ ($\Id^A_j$ if we 
need to mention explicitly the algebra $A$ in order to avoid confusion). 
We explore the question of the Cohen-Macaulay property and the dimensions
of modules $\Fl M$ and of their annihilators.

The answer to this question happens to be connected with the class of 
maximally central algebras introduced by Azumaya in \cite{AzN, Az}:

{\bf Definition } \cite[\S2]{AzN}. A finite-dimensional associative
algebra $A$ with identity over a field $\Bbbk$ is called \emph{maximally
central}, if $A$ is a direct sum of algebras $A_i,$ whose quotients modulo 
the radical are simple and
$$
\dim_\Bbbk A_i\leq t_i^2\dim_\Bbbk Z(A_i),
$$
$t_i^2$ being the rank of $A_i/\rad A_i$ over its center and actually an
equality takes place.

If $t_i$ are the same for all the summands, we call a maximally central
algebra \emph{equidimensional}.

Further equivalent definitions of these algebras are collected in
Section~\ref{defbigc}. See, \eg, \cite{DeM} for the results of further
development of the works of Azumaya that lead to the concept of Azumaya
algebras.

In the present paper we prove

\begin{stat} \label{stat1}
Suppose that either $l=1$ or $A$ is maximally central. Then
\begin{enumerate}
\item\label{two}
$\Fl{\cdot}$ is an exact fully faithful functor from the category of
fi\-ni\-te-di\-men\-sional $A$-modules to the category of graded 
$R$-modules and homogeneous homomorphisms of degree $0$;
\item\label{good} if $l=1$ or $A$ is equidimensional with $t_i=n$, then
$\Fl{\cdot}$ transforms finite-dimensional $A$-modules into Cohen-Macaulay 
$R$-modules of projective dimension  $(l-1)n+1$ (which equals 1 for 
$l=1$), and for any maximally central algebra $A$ $\Fl{\cdot}$ transforms 
indecomposable finite-dimensional $A$-modules into Cohen-Macaulay  
$R$-modules;
\item\label{three}
for $l=1$ and every $M$ the annihilator of $\Fl M$ is a principal ideal.
\end{enumerate}
\end{stat}

There is a converse statement:

\begin{stat}\label{stat2}
If for some $l>1$ either $F_l$ is exact or for all indecompsable modules 
$M$ the modules $F_l(M)$ are Cohen-Macaulay, then $A$ is a maximally 
central algebra, and if the indecomposability condition is omitted, then
$A$ is an equidimensional maximally central algebra. Furthermore, for any 
(associative unitary) algebra $A$, some $A$-module $M$ and some $l>1$ 
$F_l(M)$ is Cohen-Macaulay, if and only if $A/\ann M$ is an 
equidimensional maximally central algebra.
\end{stat}

The proof of Theorem~\ref{stat1} (namely, Lemma~\ref{Mn}) supplies
information on the minimal resolution of $\Fl M,$ which allows to determine
various invariants of these modules, in particular, to prove again and
generalize the main results of~\cite{AL}:

\begin{stat}
Suppose that a maximally central algebra $A$ is equidimensional and all 
$t_i$ equal $n.$ Then for every $A$-module $M$ the invariants of $\Fl M$ are
equal to $\dim_{\K}M/n$ times the invariants of $\Fl P$ for a simple
$\B A$-module $P,$ where $\B A=A\otimes_{\K}\B\K$ is a
$\B\K$-algebra and $\B\K$ is the algebraic closure of $\K.$
Here the invariants of $\Fl P$ have the following values:
\begin{itemize}
\item the Betti numbers $b_0=n, b_1=ln, b_i=\binom{ln}{n+i-1}\binom{n+i-3}{n-1}$
for $i\geq 2,$ concentrated in degree $0,1, n+i-1$ respectively (\ie the 
ranks of the modules $F_i$ in the minimal graded free resolution equal 
$b_i$ and each $F_i$ has generators only in one degree, namely, in degree
$0$ for $i=0,$ $1$ for $i=1,$ $n+i-1$ for $i\geq 2$);
\item the Cohen-Macaulay type $t=b_{(l-1)n+1}=\binom{ln-2}{n-1}$;
\item the Hilbert series ($\sum_k\dim_\K F_l(P)_k t^k$, $F_l(P)_k$ being 
the homogeneous component of $F_l(P)$ of degree $k$) 
$$
\Fl{P}(t)=
(1-t)^{(l-1)n+1-\dim R}\sum_{i=0}^{n-1}\binom{ln}{i}(n-i)t^i(1-t)^{n-1-i};
$$
\item the multiplicity
$$
e=(\dim R-(l-1)n-2)!\lim_{k\to\infty}\dim_\K 
F_l(P)_k/k^{\dim R-(l-1)n-2} =\binom{ln}{n-1}.
$$
\end{itemize}  
\end{stat}

Let us remark that the construction descirbed in the beginning of the
introduction makes it possible to associate some systems of PDEs with
constant coefficients to the modules over polynomial rings that we
consider, so it is possible to restate the results in terms of PDEs.
Besides, the results of this work may be applied to quaternionic analysis:
it would be interesting to obtain from the Eagon-Northcott complex, which
serves as the minimal resolution for the module $\M_n$ of~\cite{palam}, an
explicit acyclic resolution for the sheaf of quaternionic-differentiable
functions (a simile of the Dolbeault complex for holomorphic functions) and
to supply by means of this resolution Theorem~3.1 and Corollary~3.4
of~\cite{palam} with analytic proofs, that the authors lacked. Let us also
remark that in \cite{Dirac} the authors have managed to study the
corresponding modules by means of Gr\"obner bases only with machine
computations for two and three ``variables'', so maybe applying to this
problem the methods of commutative algebra similar to those used in the
present paper will result in further progress.

Chapter~\ref{prelim} is devoted to preliminaries: Sections 0.1--\ref{inv} 
deal with commutative algebra and sections \ref{decomp}--\ref{defbigc} 
with finite-dimensional associative algebras. In Section \ref{defbigc} 
equivalence of several definitions of a maximally central algebra is shown 
(Proposition~\ref{bigcenter}), of which Definitions \ref{:bigc}) and
\ref{:lincomb}) seem to be new. There is also an example to show that 
Definition \ref{:lincomb}) is equivalent to the others only in case of a 
perfect field, and otherwise it is more restrictive. Three following 
chapters contain proofs of the respectively numbered theorems.

This is an English translation of the author's Ph.\:D. thesis defended in
December 2004 at The Department of Mechanics and Mathematics of the
Moscow State (Lomonosov) University. The results have been published in
\cite{art0}, \cite{art1}, \cite{art2}, \cite{art3}, \cite{art4kurz}, but the 
proof of the second part of Theorem~\ref{stat2} (the one concerning just 
one module) appears here for the first time. Also the exposition here is 
more compact and straightforward as in the papers above, unnecessary
repetitions have been removed. I would like to thank my mentor 
Prof.~E.~S.~Golod for his constant attention and instructive remarks, the 
anonymous referee of~\cite{art1} for simplifying some arguments and
also A.~A.~Gerko for a motivation to finish this work.

\chapter*{Notations}
\addcontentsline{toc}{chapter}{Notations}

\hspace{\parindent}$\set R$ --- the field of real numbers,

$\set C$ --- the field of complex numbers,

$\set H$ --- the skew-field of quaternions,

$\BF_p$ --- the finite field with $p$ elements,

$M_n(K)$ --- the ring of $n\times n$ matrices with elements in a ring $K,$

$\K$ --- a field,

$A$ --- a finite-dimensional associative algebra with identity over a field
$\K,$

$\B\K$ --- the algebraic closure of $\K,$

$\B A=A\otimes_\K\B\K,$

$f_1,\dots,f_d$ --- a vector space basis for $A$ over $\K,$

$R=\K[x_{11},\dots,x_{dl}]$ --- the ring of polynomial functions on an 
affine space $A^l,$

$\rho_M$ --- the matrix representation of an algebra $A$ corresponding to a 
finite-dimensional $A$-module $M,$

$\Id_j$ or $\Id^A_j$ --- a generic element of an algebra $A$, see 
Introduction for a coordinate description, the beginning of Chapter 1 for 
an invariant one,

$A|B$ --- an $(m+n)\times k$ matrix obtained by writing the $n\times k$ 
matrix $B$ to the right of the $m\times k$ matrix $A$,

$\Fl\cdot$ or $F_l^R(\cdot)$ or $F_l^K(\cdot)$ --- the functor we study,
see Introduction for a coordinate description, the beginning of Chapter 1 for 
an invariant one; here $R$ (maybe with primes) denotes the polynomial ring 
used in the construction, $K$ (some letter that is not $R$ with primes) 
the finite dimensional algebra over which the construction is performed,

$S^i G$ --- a symmetric power (of a free module over a commutative ring),

$S(G)$ --- the symmetric algebra of a vector space,

$\bigwedge^i G$ --- an exterior power (of a free module over a commutative 
ring), 

$G^*$ --- the dual module (of a free module over a commutative ring),

$\dim M$ --- the Krull dimension (of a module $M$ over a polynomial ring),

$\mathop{\rm supp} M$ --- the support of a module over a commutative ring,

$\dim_\K M$ --- the dimension of a vector space $M$ over $\K,$

$l(M)$ --- the length of a module $M$ over a finite-dimensional algebra,

$\depth M$ --- the depth of the homogeneous maximal indeal in a polynomial 
ring on a graded module over this ring (see Section~\ref{dpcm}),

$Q(K)$ --- the field of fractions of a commutative ring $K$,

$k(\fp)=Q(K/\fp)$ --- the residue field of a prime $\fp$ in a commutative 
ring $K$,

$\pd M$ --- the projective dimension of a module over a polynomial ring,

$M_i$ --- the homogeneous component of degree $i$ of a graded module $M$ 
over a polynomial ring,

$M[i]$ --- a grading shift for a graded module over a polynomial ring 
($M[i]_j=M_{i+j}$), 

$M(t)\in\set Z[[t]]$ --- the Hilbert series of a graded module $M$ over a 
polynomial ring,

$\hht I$ --- the height of an ideal in a commutative ring,

$\Ass M$ --- the set of associated primes for a module over a commutative 
ring,

$\ann M$ --- the annihilator of a module,

$Z(K)$ --- the center of a ring $K$,

$K^0$ --- the opposite ring of a ring $K$ (the additive group is the same 
as in $K$, and the product $ab$ in $K^0$ equals the product $ba$ in $K$),

$\rad K$ --- the radical of a finite-dimensional algebra $K,$

$\Br(F)$ --- the Brauer group of a field $F$,

$\Br(F, L)$ --- the subgroup in $\Br(F)$ comprising the classes of central 
simple algebras that split over the extension $L$ of $F$,

$H^2(G,K^*)$ --- the second group cohomology of a group $G$ with 
coefficients in the multiplicative group of a field $K.$

\chapter{Preliminaries}\label{prelim}

\section{Depth and Cohen-Macaulayness}\label{dpcm}
We use here the graded versions of these concepts, analogous to the local
ones considered in \cite{Se}, where depth is called ``codimension
homologique". An exposition of the graded case can be found in \cite[\S 
1.5]{BH}, but we tried to give references to sources available in Russian 
wherever possible.

\mydef\label{regseq}(\cite[Chap.~17, p.~423]{Ei}; \cite[chap.~IV, A.4]{Se},
\emph{$M$-suite.}) Let $R$ be a commutative Noetherian ring, $M$ an
$R$-module. A sequence $a_1,\dots,a_n\in R$ is called \emph{$M$-regular,} if
$(a_1,\dots,a_n)M\ne M$ and for $i$ between 1 and $n$ the multiplication
by $a_i$ on the module $M/(a_1,\dots,a_{i-1})M$ is injective.

\begin{prop} {\rm\cite[Theorem~17.4]{Ei}}
Let $I$ be an ideal in $R$ with $IM\ne M.$ Then all maximal $M$-regular
sequences in $I$ (\ie those that cannot be continued) have the same length.
\end{prop}

\mydef(\cite[17.2, p.~429]{Ei}, \cite[chap.~IV, A.4, d\'efinition~6]{Se}.)%
\label{dp}
Under the hypotheses of the previous proposition the length $\depth(I,M)$
of any maximal $M$-regular sequence is called \emph{the depth} of $I$ on
$M.$

When considering the depth of the ideal of the polynomials vanishing at the
origin on graded modules over a polynomial ring, we shall omit the ideal
and talk of the depth of a module.

\begin{prop}\label{dpft}{\rm(\cite[Theorem~16.11]{BV}, \cite[Prop.~18.4,
Cor.~18.6]{Ei}, \cite[chap.~IV, A.4, prop.~6~ff.]{Se}.)}
For every finitely generated $R$-module $M$ one has:

$\depth(I,M)\leq \dim M,$ $\dim M$ being the Krull dimension;
$$\depth(I,M)=\min\{i\mid\Ext^i_R(N,M)\ne 0\}$$
for a finitely generated module $N$ with the support equal to the closed
subset defined by $I,$ so that the depth does not change if one replaces
$I$ by its radical;
$$\depth(I,M\oplus N)=\min\{\depth(I,M),\depth(I,N)\};$$
in a short exact sequence $0\to M\to N\to P\to 0$ one has
\begin{align*}
\depth(I,N)&\geq\min\{\depth(I,M),\depth(I,P)\},\\
\depth(I,M)&\geq\min\{\depth(I,N),\depth(I,P)+1\}.
\end{align*}
\end{prop}

\begin{prop}[the Auslander-Buchsbaum formula]\label{dpft2}
{\rm(\cite[chap.~IV, D.1, prop.~21]{Se}, \cite[Ex\-er\-cise~19.8]{Ei}.)}
For a graded module $M$ over a polynomial ring $R$ one has
$\depth M+\pd M=\dim R.$
\end{prop}

\mydef\label{CM}(\cite[chap.~IV, B.1, d\'ef.~1]{Se},
 see also \cite[Chap.~18, p.~451]{Ei}.)
A module $M$ over $R$ is called \emph{Cohen-Macaulay,} if for every
maximal ideal $\mm$ in $R$ one has $\depth(\mm,M)=\dim M.$

\begin{prop}\label{grad.CM}{\rm\cite[Prop.~16.20]{BV}}
Let $M$ ba a graded module over a polynomial ring $R.$ Then $M$ is
Cohen-Macaulay iff $\depth M=\dim M$ (the depth of the homogeneous maximal
ideal being considered).
\end{prop}
\begin{prop}\label{regCM}{\rm\cite[chap.~IV, B.2]{Se}}
Let $M$ be Cohen-Macaulay. Then a sequence $a_1,\dots,a_s$ is $M$-regular iff
$M/(a_1,\dots,a_s)M\ne 0$ and
$$
\dim M/(a_1,\dots,a_s)M =\dim M-s,
$$
(where for a Cohen-Macaluay $M$ and any sequence of length $s$ factorizing 
$M$ modulo the sequence reduces its dimension for at most $s$), and then 
the quotient module $M/(a_1,\dots,a_s)M$ is Cohen-Macaulay.
\end{prop}

\section{The Eagon-Northcott Complex}\label{WMCE}
\mydef \label{EN}(\cite[2.C]{BV}, $\mathcal{D}_1(g)$; \cite[A2.6.1, p.~600]
{Ei}, $\CC^1,$ \emph{the Buchsbaum-Rim complex.})
Let $\varphi=(\varphi_{ij})$ be a $(g\times f)$-matrix determining a
homomorphism of free modules $\varphi\co F\to G,$ and let the ranks of these
modules be denoted by the corresponding lowercase letters. \emph{The number 1 
Eagon-Northcott complex} constructed from $\varphi$ is the complex
$$
0\to F_k\myto{\varphi_k} \dots \myto{\varphi_2} F_1\myto{\varphi_1}F_0
$$
of free modules, where $k=f-g+1,$ $F_0=G,\ F_1=F, 
\ F_i=\bigwedge^{g+i-1}(F)\otimes S^{i-2}(G)^*$ for $i\geq 2,$ and the
differentials are as follows. Let $f_1,\dots,f_f$ and $g_1,\dots,g_g$ be
the bases of $F$ and $G$ respectively, then
$$
\varphi_1(f_j)=\sum\varphi_{ij}g_i,\quad
\varphi_2(f_{j_0}\wedge\dots\wedge f_{j_g})=
\sum_{k=0}^g (-1)^k M_{j_0\dots\widehat{j_k}\dots j_g}f_{j_k}
$$
(here and in the display below $\widehat\cdot$ denotes omitting an item
from a list), $M_{j_1\dots j_g}$ being the $(g\times g)$-minor of $\varphi$
consisting of the columns $j_1,\dots,j_g$ in the specified order, and
$$
\varphi_i(f_{j_0}\wedge\dots\wedge f_{j_{g+i-2}}
\otimes (g_{i_1}\dots g_{i_{i-2}})^*)=$$$$
=\sum_{k=0}^{g+i-2}(-1)^k\sum_{k'=1}^{i-2}
\varphi_{i_{k'}j_k}f_{j_0}\wedge\cdots\widehat{f_{j_k}}\dots\wedge 
f_{j_{g+i-2}}
\otimes (g_{i_1}\dots\widehat{g_{i_{k'}}}\dots g_{i_{i-2}})^*.$$

\remm In \cite{Ei} similar number $l$ complexes for all integers $l$ are
defined, but we need only this particular case.

\medskip
For such complexes one has the following exactness criterion.

\begin{prop}\label{prop.WMCE}
{\rm(\cite[Theorem~16.15]{BV}, see also \cite[Theorem]{WMCE},
\cite[Chap. 20.3]{Ei}.)}
Let $R$ be a commutative Noetherian ring, $M\ne 0$ a finitely generated
$R$-module,
$$\mathbb{A}=(0\to
F_k\myto{\varphi_k}F_{k-1}\to\ldots\to F_1\myto{\varphi_1}F_0)$$
a complex of finitely generated free $R$-modules. Let
$$r(j):=\sum_{i=j}^k(-1)^{i-j}\rk F_i$$ be the rank (\ie the order of
the maximal nonsero minor) of $\varphi_j$ in the case when this complex is
exact. Let $I_r(\varphi)$ denote the ideal generated by all the 
$(r\times r)$-minors of $\varphi.$
Then $\mathbb{A}\otimes_R M$ is exact, iff for all $j$ between 1 and $k$
$I_{r(j)}(\varphi_j)$ contains an $M$-sequence of length $j$ or
$I_{r(j)}(\varphi_j)M=M.$
\end{prop}
 
The following fact offers a great simplification of this criterion for the
Eagon-Northcott complex:
\begin{prop}{\rm\cite[Theorem A2.10 b]{Ei}}\label{prop2.WMCE}
For any matrix $\varphi$ the rank of $\varphi_j$ in the corresponding
Eagon-Northcott complex does not exceed $r(j),$ and the ideal 
$I_{r(j)}\varphi_j$ of the rank $r(j)$ minors of 
$\varphi_j$ lies in the ideal $I_m(\varphi)$ of the maximal minors of 
$\varphi$ and has the latter for its radical.
\end{prop}

Therefore, to prove the exactness of a complex of type
$\mathbb{A}\otimes_R
M,$  $\mathbb{A}$ being the Eagon-Northcott complex, it suffices to check
that $I_m(\varphi)$ contains an $M$-sequence of length $f-g+1$: then 
Propositions~\ref{prop2.WMCE} and~\ref{dpft} show this to be the 
case for all $I_{r(j)}(\varphi_j),$ so Proposition~\ref{prop.WMCE} applies.

\section{The Invariants of Modules over Commutative Rings}
\label{inv}
\mydef (cf. \cite[Exercise~A3.18]{Ei})
Let $R=\bigoplus_{i\geq 0}R_i$ be a graded Noetherian ring, $R_0=\K$ a field,
 $R_+=\bigoplus_{i>0}R_i$ the maximal homogeneous ideal,
$R[n]$ the free $R$-module with shifted grading (the generator having degree
$-n$). Let $M$ be a finitely generated graded $R$-module. It is known
 \cite[19.1; Theorem~20.2]{Ei}, that in this case there exists a unique free
resolution 
$$
\dots\to F_i\to\dots\to F_1\to F_0\to M\to 0
$$
over $R,$ such that
$F_i=\bigoplus_jR[-j]^{b_{ij}},$ the homomorphisms are homogeneous of degree
0 and $d(F_i)\subset R_+F_{i-1}.$
\emph{Graded Betti numbers} of a graded module $M$ are the
 $b_{ij},$ \emph{Betti numbers} $b_i=\sum_jb_{ij}=\dim_{\K}\Tor_i^R(\K,M).$ 
If for some $i$ $b_{ij}=0$ for every $j$ except $j_i,$ we say that the $i$th
Betti number is \emph{concentrated in degree} $j_i.$

\mydef\cite[Exercises~10.12--10.13]{Ei}
\emph{The Hilbert series} $M(t)$ of a graded module $M$ is the series
$\sum\dim_{\K}M_it^i.$ 

Then one has $M(t)=p(t)/(1-t)^{\dim M}$ for a finitely generated module $M$
over a polynomial ring and a polynomial $p(t)$ with $p(1)\ne 0.$ The
Hilbert series is additive in short exact sequences of graded modules and
homogeneous homomorphisms of degree $0,$ which follows from the additivity
of dimension for homogeneous components of each degree. In particular, if
$x$ is a homogeneous nonzerodivisor of degree $d$ on a graded module $M$,
the exact sequence
$$
0\to M[-d]\xto{x} M\to M/xM\to 0
$$ 
shows that $(M/xM)(t)=(1-t^d)M(t).$ Induction derives from this a formula for 
the Hilbert series of the quotient module modulo a homogeneous regular 
sequence, which we shall use.

\mydef \emph{The multiplicity} of a module $M$ (\wrt $R_+$) is $e=p(1).$
The definition shows the mutiplicity to be also additive in short exact
sequences. One can easily check that for graded modules this number
coincides with the one introduced in \cite[chap.~V, A.2]{Se} and
\cite[12.1]{Ei} (the last definition was included in the statement of 
Theorem~3). 

\mydef (cf. \cite[Exercise~21.14]{Ei}) 
\emph{The Cohen-Macaulay type} of a Cohen-Macaulay module $M$ is the number 
 $t(M)=\dim_{\K}\Ext^{\depth M}_R(\K,M),$
where $\K=R/R_+$ as an $R$-module.
\begin{prop}\label{CMt}
Over a polynomial ring $R$ one has $t(M)=b_{\pd_RM}.$ 
\end{prop}
\Proof Over a regular ring $\Tor^R_i(\K,M)$ is isomorphic to 
$\Ext^{\dim R-i}_R(\K,M)$ 
\cite[chap.~IV, D.1, cor. 1 au th\'eor\`eme 5]{Se}, and the
Auslander-Buchsbaum formula yields the required equality. \qed

\section{Gr\"obner Bases}

As most expositions of the theory of Gr\"obner bases treat the case of ideals,
not submodules as used here, we recall the basic statements from
\cite[Chap.~15]{Ei}.

Let $R$ be a polynomial ring over a field $\K,$ and $F$ a free $R$-module 
with a chosen basis $e_1,\dots,e_s.$ A \emph{monomial} in $F$ is an element 
of the form $me_i,$ $m$ being a monomial in $R$ (\ie a product of powers 
of the variables). A \emph{monomial order} in $F$ is a total order on the
set of all monomials in $F$ such that if $m_1>m_2$ are two monomials in $F$
and  $n\ne 1$ is a monomial in $S$, one has $nm_1>nm_2>m_2.$ Every such
order is Artinian (every nonempty subset has a least element).

Let us describe several ways to construct these orders that we shall use.
Take any total order on the set of variables in $R$. Then one can induce
the lexicographic order on monomials in $R$: compare the degrees \wrt the
greatest variable, if they are equal, proceed to the next variable, and so
on. One can induce the degree-lexicographic order: first compare the total 
degree of monomials, and in case of equality compare them lexicogaphically.
The same can be done for noncommutative polynomials, where monomials are 
words, so the lexicographic order compares them letter-by-letter (in the
noncommutative case, not all orders are Artinian, but this one is). Given 
an order on the monomials in the ring and a total order on the basis 
elements of a module, there are two ways to construct an order on the 
monomials in the module: ``term over position'', comparing first the 
coefficients \wrt the order in the ring and then basis elements in case of 
equality, and ``position over term'', comparing first the basis elements 
and then the coefficients.

\emph{The initial term} of an element $f=\sum a_im_i\in F,$ for 
$a_i\in\K^*$ and $m_i$ being different monomials in $F$, is $a_0m_0,$ $m_0$ 
being the greatest of the $m_i$ involved. If $M\subset F$ is a submodule, 
the \emph{initial module} of $M$ is the submodule in $F$ generated by the 
initial terms of all the elements in $M.$ Then the images of the monomials in 
$F$ not contained in the initial submodule of $M$ constitute a vector 
space basis for $F/M$. In particular, a homogeneous $M$ has the same 
Hilbert series as its initial submodule, so the same is true for the 
quotient modules \cite[Theorem~15.26]{Ei}. 

Suppose a monomial $n$ involved in an element $g$ with coefficient $a$ is 
divisible by the initial term $m$ of a element $f$ (multiplying $f$ by a 
scalar one can assume that the coefficient at $m$ in $f$ equals $1$). Then
the \emph{reduction} of $g$ by $f$ is replacing $g$ with $g-a(n/m)f.$ If 
no monomial in $g$ is divisible by the initial terms of $f_1,\dots,f_k,$ 
$g$ is said \emph{not to be reducible} by $f_1,\dots, f_k.$ The remark in 
Section~\ref{defbigc} uses a version of reduction for two-sided ideals in 
a noncommutative polynomial ring: divisibility means that $n=n_1mn_2$ and 
reduction replaces $g$ with $g-an_1fn_2.$ The Artinian property of the 
order ensures that there can be no infinite sequence of reductions.

A set of elements $g_1,\dots,g_k$ in module $M$ is called a
\emph{Gr\"obner basis} for $M$ \wrt a given order, if the initial terms of 
these elements generate the initial submodule of $M$ (in the homogeneous 
case this condition can be verified by means of Hilbert functions). Then 
these elements generate $M.$

If $m_1$ and $m_2$ are two monomials in $F$ containing the same basis
element $e_i,$ the least common multiple $m$ of these monomials is defined
in an obvious fashion. If these monomials are the initial terms of $f,g\in
F$ respectively, the \emph{S-polynomial} corresponding to the pair $(f,g)$
is the element $(m/m_1)f-(m/m_2)g\in F.$

Buchberger's criterion \cite[Theorem~15.8]{Ei} says that the set of 
elements $g_1,\dots,g_k\in F$ is a Gr\"obner basis of the submodule they 
generate iff any S-polynomial corresponding to a pair of elements in this 
set (with the same basis element of $F$ in their initial terms) can be 
reduced to zero by a sequence of reductions by the elements of this set, 
and then if we apply arbitrary reductions by $g_i$ to an S-polynomial and 
after some steps obtain a element not reducible by $g_1,\dots,g_k$, this 
element is zero.

\section{The Direct Sum Decomposition of Fi\-ni\-te-Di\-men\-sional
Algebras}\label{decomp}
A finite-dimensional associative algebra $A,$ considered as a left module
over itself, can be decomposed into a direct sum of indecomposable left
ideals \cite[Theorem~14.2]{CR}, and the decomposition into a direct sum of
subalgebras is obtained from this one by grouping the summands
 \cite[\S\S54--55]{CR}. Namely, two such ideals $a$ and $b$ are called
\emph{linked} \cite[Def.~55.1]{CR}, if there is such a chain of
indecomposable left ideals $a=a_1, a_2, \dots, a_n=b,$ that any two
neighbouring ideals have a composition factor in common. Direct sums of
linkage classes --- called blocks --- are indecomposable two-sided ideals in
$A,$ they are uniquely determined and $A$ is the direct sum of them
 \cite[Theorem~55.2]{CR}.

We shall apply this result in the case when every indecomposable $A$-module
has only 1 type of composition factors. Then every block has only 1 type
of composition factors, so the algebra is a direct sum of algebras, where
each summand has only 1 simple module.

Let us as well recall the description of direct sum decompositions of an 
algebra in terms of idempotents \cite[\S 25, exercise~2]{CR}: there is a 
one-to-one correspondence between the decompositions of algebra into a 
direct sum of two-sided ideals and the decompositions of identity into a 
sum of central orthogonal idempotents, \ie idempotents from the center of 
the algebra with all pairwise products equal to $0.$ The correspondence is 
natural: idempotents are taken to the ideals they generate and the 
idempotents corresponding to a direct sum decomposition are the 
projections of identity into the summands.
 
\section{Simple and Separable Algebras}

We write $Z(B)$ for the center of a finite-dimensional algerba $B$, $\rad
B$ for its radical, $B^0$ for the opposite algebra (\ie the same vector
space with the multiplication $(b_1,b_2)\mapsto b_2b_1$).

\begin{prop}{\rm\cite[Theorem~68.1]{CR}}\label{Br}
Let $S$ be a simple algebra, $K$ its center and $L$ a field extension of $K.$ 
Then $S\otimes_KL$ is a simple algebra with center $L.$
\end{prop}

\mydef \cite[Def.~71.1]{CR} A semisimple finite-simensional algebra over a
field is called \emph{separable,} if it remains semisimple after any
extension of the base field. 

\remm\label{69.8} Any semisimple algebra over a perfect field is separable,
because Theorem~69.4 of \cite{CR} guarantees it to remain semisimple after
any finite (separable) extension of the base field, hence after taking the
algebraic closure, for otherwise the nilradical will be defined over some
finite extension. So over the algebraic closure of the base field this
algebra is isomorphic to a direct sum of matrix ones, and Theorem~71.2 of
\cite{CR} says that an algebra isomorphic over some extension of the base
field to a direct sum of matrix algebras is separable.

\begin{prop}[Wedderburn--Malcev Theorem]{\rm\cite[Theorem~72.19]{CR}}%
\label{WeddMal}
Let $B$ be a finite-dimensional algebra over a field, such that $B/\rad B$
is a separable algebra. Then $B$ contains a subalgebra $S,$ so that one has 
a semi-direct sum decomposition $B=S\oplus\rad B.$
\end{prop}

\begin{prop}{\rm\cite[Chap.~IV, Theorem~4.4.2]{Her}}\label{factor}
If $S\subset A$ is a finite-dimensional simple subalgebra of an 
$L$-algebra $A$ containing the identity element of $A$ and $Z(S)=L,$ then
$A=S\otimes_{Z(S)}K,$ $K$ being the centralizer of $S$ in $A.$
\end{prop}

The following well-known lemmas will be used in the sequel.
\begin{lm}\label{dimc}
Suppose that all the simple quotients of an algebra $A$ are of dimension
$n^2$ over their centers. Then the same is true for the algebra $\B A$ 
obtained from $A$ by extending the base field to its algebraic closure.
\end{lm}
\Pf 
As an extension of the base field takes nilradical into nilradical, one 
can suppose that $A$ is semisimple, and considering every summand in turn 
one can suppose that is is simple. Let $\B\K$ be the algebraic closure of 
$\K.$ Then
$$
A\otimes_\K\B\K/\rad\B A=A\otimes_{Z(A)}((Z(A)\otimes_\K\B\K)/\rad 
(Z(A)\otimes_\K\B\K)),
$$
as the tensor product of a central simple algebra and a semisimple one is 
semisimple \cite[Theorems 68.1 and 71.10]{CR}. But the second factor in 
this formula is a commutative semisimple algebra over $\B\K,$ \ie the sum 
of several copies of $\B\K,$ and $A\otimes_{Z(A)}\B\K$ is a simple algebra 
of the same dimension over its center \cite[Theorem~68.1]{CR}. \qed

\begin{lm} \label{matr}
Let $A$ be a finite-dimensional associative unitary algebra over an 
algebraically closed field $\K$ with its simple quotients of dimension 
$n^2$ over their centers. Then $A=M_n(K)$ for some algebra $K.$
\end{lm}
\Pf
As the field is algebraically closed, the semisimple quotient of $A$ is
the sum of several copies of $M_n(\K)$. Proposition~\ref{WeddMal} shows
that this quotient can be embedded as a subalgebra of $A$, so that $A$ is a
semi-direct sum of this subalgebra and the radical of $A$. As the radical
is nilpotent, it is easy to see that the subalgebra contains the identity
element of the whole algebra. If we embed $M_n(\K)$ as the diagonal of this
subalgebra, then the matrix algebra will be also embedded in $A$ as a
subalgebra containing the identity of $A$. As the center of $M_n(\K)$
equals $\K,$ one can apply Proposition~\ref{factor} and conclude that
$A=M_n(\K)\otimes_\K K=M_n(K).$ \qed

\section{Maximally Central Algebras}
\label{defbigc}
{\bf Definition} \cite[\S2]{AzN}. A finite-dimensional associative  
algebra $A$ over $\Bbbk$ with identity element is called \emph{maximally 
central}, if $A$ is a direct sum of algebras $A_i$ the semisimple 
quotients of which are simple and
$$
\dim_\Bbbk A_i\leq t_i^2\dim_\Bbbk Z(A_i),
$$
where $t_i^2$ is the rank of $A_i/\rad A_i$ over its center and actually 
an equality takes place.

If the $t_i$ are the same for all the summands, we call a maximally 
central algebra \emph{equidimensional}.

The definition of maximally central algebras is generalized in
 \cite{Az} to algebras over a Henselian ring that are free modules of 
finite rank over that ring, but we restrict ourselves to the case of 
finite-dimensional algebras over a field. This definition allows for a lot 
of equivalent restatements (in particular, the main results of this paper 
can be interpreted as such a restatement), mostly given in \cite{AzN,Az}:

\begin{prop}\label{bigcenter}
The following conditions on a finite-dimensional associative algebra $A$ 
with $1$ over a field $\Bbbk$ are equivalent:
{
\let\wwww\labelenumi
\renewcommand{\labelenumi}{\textup{\wwww}}
\begin{enumerate}
\item \label{:maxc} 
$A$ is a maximally central algebra;
\item \label{:bigc} 
if $A$ is projected onto any of its quotient algebras $B$, the center of $A$ 
is mapped surjectively onto the center of $B$;
\item \label{:barbigc} 
the algebra $\B A$ obtained from $A$ by extending the base field to its 
algebraic closure is a direct sum of matrix algebras over their centers;
\item \label{:extmaxc}
the algebra $A\otimes_\Bbbk L$ is maximally central, $L$ being an arbitrary
field extension of $\Bbbk$;
\item \label{:Azum1}
$A$ is a direct sum of algebras $A_i$ with centers $Z_i$ in such a way 
that $A_i$ are free $Z_i$-modules and 
$\End_{Z_i}A_i=A_i\otimes_{Z_i}A_i^0$;
\item \label{:Azum}
$A$ is an Azumaya algebra over its center $Z,$  \ie 
\textup{(\cite[Chap.~VI, \S41, Definition 41.5]{Ma})} a projective
$Z$-module such that for any prime ideal $\fp\subset Z$
$A\otimes_Zk(\fp)$ is a central simple algebra over $k(\fp).$

If the base field is perfect, these condintions are equivalent to the 
following:

\item \label{:lincomb}
$A=\bigoplus_i S_i\otimes_{Z(S_i)}K_i,$ where $S_i$ are simple algebras 
over $\K$ and $K_i$ are commutative algebras over the corresponding
$Z(S_i)$ ($A$ is ``a linear combination of simple algebras with 
commutative coefficients'').
%
% Видимо, переводы из Серра вредно отражаются на стиле.
%
\end{enumerate}
}
\end{prop}

\cite[Chap.~2, \S\S 2,3]{DeM} contains some equivalent definitions of
Azumaya algebras, allowing to give some more restatements of
condition~\ref{:Azum}).

\Pf The equivalence of \ref{:maxc}), \ref{:barbigc}) and \ref{:extmaxc})
is proved in \cite[\S2, Theorem 2, Corollary]{AzN}. The equivalence of 
\ref{:maxc}) and \ref{:Azum1}) is proved in \cite{Az} at the end of Section~4.

The equivalence of \ref{:Azum1}) and \ref{:Azum}) follows from Theorem~15 
of \cite{Az}. This theorem states that a $Z_i$-algebra $A_i$ that is a 
free $Z_i$-module, $Z_i$ being a commutative ring,
satisfies the property required from the algebra denoted this way in 
condition~\ref{:Azum1}), iff for every maximal ideal $\fp\subset Z_i$ 
$A_i/\fp A_i$ is a central simple algebra over $Z_i/\fp.$ The center of 
$A$ is a commutative Artinian ring that can be decomposed into a sum of 
local Artinian rings, and then the same idempotents give a decomposition 
of $A$ into a direct sum of algebras with local centers, so the summands 
are indecomposable. The condition \ref{:Azum}) is local \wrt the center, 
so it can be transferred to every summand and, as a finitely generated
projective module over a local ring is free, it means that the summands are
free over their centers. The condition~\ref{:Azum1}) can also be
transferred to summands, if one decomposes both sides of the equation into
a sum of the modules over the summands of the center. Now the equivalence of
these conditions for one summand is claimed in Theorem~15 of \cite{Az}, as
quoted above.

\ref{:barbigc})\hence\ref{:bigc}): given a quotient algebra of $A$, it is 
enough to check the surjectivity of the map of the centers for the 
correponding quotient of $\B A$ after passing to the algebraic closure of 
the field: (a faithfully flat) base change commutes with taking the center,
as the latter is defined by linear equations. So it suffices to establish 
this condition for $\B A$. This algebra is a direct sum of matrix algebras 
with commitative coefficients, so any two-sided ideal in it is a direct 
sum of ideals in the summands, and any two-sided ideal in a matrix algebra 
over a ring is the set of matrices with elements from some two-sided ideal 
in the ring (see \cite[chap.~VIII, \S7, exercice 6 b)]{Bour}, where it is 
presented in generality proper to the author). So any quotient algebra is 
a direct sum of matrix algebras over some quotients of their centers. As 
the center of a matrix algebra over a commutative ring coincides with that 
ring, the map of the centers is surjective.

\ref{:bigc})\hence\ref{:maxc}): let us decompose the semisimple quotient 
of $A$ as a direct sum of simple algebras. This decomposition is given by 
a complete family of central orthogonal idempotents in the semisimple 
quotient. As the center of $A$ maps onto the center of the semisimple 
quotient by hypothesis, Theorem~24 of \cite{Az} (or Corollary~7.5 of
\cite{Ei}) allows one to lift this family to a complete family of 
orthogonal idempotents in the center of $A$, the latter giving a 
decomposition of $A$ into a direct sum of algebras, each of which has only 
one simple quotient and satisfies the condition \ref{:bigc}).

For every summand $A_i$ we show the formula by induction on Loewy length, 
\ie the nilpotence degree of the nilradical. If the nilradical is zero,
$A_i$ is simple and its center is a field of, obviously, required dimension.
Otherwise let $(\rad A_i)^n$ be the last non-zero power of the radical. 
Then by induction assumption the center of $A_i/(\rad A_i)^n$ has the 
dimension required. The center of $A_i$ is mapped onto it surjectively 
with kernel equal to $Z(A_i)\cap (\rad A_i)^n,$ and
$$
\dim_\K A_i-\dim_\K A_i/(\rad A_i)^n=\dim_\K (\rad A_i)^n,
$$
so it suffices to show that
$$
\dim_\K Z(A_i)\cap (\rad A_i)^n=\dim_\K (\rad A_i)^n/\dim_{Z(S_i)}S_i
$$ 
for $S_i=A_i/\rad A_i.$ The $(\rad A_i)^n$ under consideration is an
$A_i$-bimodule, \ie an $A_i\otimes_\K A_i^0$-module, but we can notice
that, as the left action of $Z(A_i)$ is the same as the right one, it is
actually a module over $A_i\otimes_{Z(A_i)}A_i^0,$ and, as both actions of
the radical are trivial and $Z(A_i)$ is mapped onto the center of the
(semi)simple quotient $S_i,$ it is actually a module over the simple
algebra $S_i\otimes_{Z(S_i)}S_i^0,$ \ie a direct sum of bimodules
isomorphic to $S_i.$ The intersection of the center of $A_i$ with the power
of the radical under consideration is the set of the elements on which two
actions of $S_i$ in this bimodule structure coincide, \ie the direct sum of
the centers of $S_i.$ Hence this intersection has the dimension required.

Now we show that these conditions are equivalent to the last one over a
perfect field. It can be easily seen that \ref{:lincomb})\hence\ref{:maxc})
over any field. Let us derive condition \ref{:lincomb}) from the other ones
over a perfect field. By \ref{:bigc}) the map from $Z(A)$ into the center
of the semisimple quotient of $A$ is surjective, so the system of central
orthogonal idempotents giving the decomposition of the semisimple quotient
into a sum of simple algebras can be lifted to a system of central
orthogonal idempotents giving a decomposition of $A$ into a direct sum of
algebras $A_i$ the semisimple quotients $S_i$ of which are simple. As the
field is perfect, by Remark~\ref{69.8} every semisimple algebra over it is
separable, so by Wedderburn--Malcev Theorem $Z(S_i)$ is embedded into
$Z(A_i)$ as a subalgebra with $1,$ so $A_i$ can be regarded as an algebra
over the field $Z(S_i).$ A finite extension of a perfect field, $Z(S_i)$ is
itself perfect, so we can apply the Wedderburn--Malcev theorem over it and
obtain an embedding of $S_i$ into $A_i$ as a $Z(S_i)$-subalgebra with~$1.$
Now we can apply Proposition~\ref{factor} and obtain that
$A_i=S_i\otimes_{Z(S_i)}K_i$, $K_i$ being the centralizer of $S_i$ in
$A_i,$  so that $K_i\supset Z(A_i).$ The tensor product decomposition
yields that $\dim_\Bbbk K_i=\dim_\Bbbk A_i/t_i^2,$ and then the inequality
from the definition of a maximally central algebra and the inclusion yield
that $Z(A_i)=K_i$, so $K_i$ is commutative. \qed

\remm\label{barequidim} Applying condition~\ref{:barbigc}) and 
Lemma~\ref{dimc} to every $A_i$ from the definition of a maximally central 
algebra, we obtain that equidimensional maximally central algebras are 
exactly those that become isomorphic to a direct sum of matrix algebras of 
the same rank over commutatvie ones after extending the base field to its 
algebraic closure, \ie just to a matrix algebra over a commutative one.

\remm Over a non-perfect field the last condition is not equivalent to the 
previous ones, as the following example shows.

Set $\K=\BF_p(x)$ for a variable $x$. Consider a purely inseparable 
extension $F=\BF_p(t)=\K[t]/(t^p-x)$ of this field. First we construct a 
finite-dimensional central skew-field over $F$ that cannot be obtained 
from such a field over $\K$ by extension of scalars.

Let $L=\BF_p(u)=F[u]/(u^p-u-t)$ be a cyclic Galois extension of degree $p$
over $F,$ $\sigma$ a generator of the Galois group taking $u$ to $u+1.$ 
Then, according to \cite[\S 114]{VDW}, the elements of the Brauer group of
$F$ that are trivial over $L$ are represented by cyclic algebras, \ie
algebras of the form $F\langle\sigma,u\rangle
/(\sigma^p-\alpha,u^p-u-t,\sigma u-(u+1)\sigma)$ (the variables commuting 
with $F$ but not with one another), $\alpha$ being an element of $F^*$ 
\cite[\S 94, 4]{VDW}. Furthermore, two such algebras are isomorphic iff 
the corresponding $\alpha$'s differ by a factor that is a norm in $L/F,$ in
particular, a cyclic algebra is the total ring of matrices iff $\alpha$ is
a norm \cite[\S 114, Aufgabe 3]{VDW}. Let us consider a cyclic algebra with
$\alpha=t-1$: it is a central simple algebra of dimension $p^2$ over $F,$
so it is either a skew-field or a matrix algebra over $F.$ As $t-1=u^p-u-1$
is an irreducible polynomial in $\BF_p[u],$ it is no norm (\ie no product
of $p$ conjugates) in the extension $L/F,$ so our cyclic algebra is a
non-trivial element of the Brauer group, so it is a skew-field. If this
skew-field could be obtained from a skew-field over $\K$ by extension of
scalars, the skew-field over $\K$ would have rank $p^2$ and would lie in
the $p$-torsion of $\Br(\K)$ by Theorem~4.4.5 of \cite{Her}, which claims
that the class of a skew-field in a Brauer group is annihilated by the
square root of the rank of this skew-field. But the map $\Br(\K)\to \Br(F)$
induced by the extension of scalars takes $p$-torsion to zero, as we are
going to show. To this end, we use the fact that the Brauer group $\Br(\K)$
is the union of its subgroups $\Br(\K,K)$ formed by algebras that become
isomorphic to a matrix algebra after tensoring with $K,$ $K$ varying over
all finite Galois extensions of $\K$ \cite[\S\S113--114]{VDW}. And
$\Br(\K,K)$ is isomorphic to $H^2(G,K^*),$ $G$ being the Galois group of
$K/\K$ \cite[\S6, n\degr8, prop.~11]{BourX} (this is essentially a
reformulation of the description of Brauer group in terms of systems of
factors). But, according to \cite[\S6, n\degr8]{BourX}, group cohomology is
an instance of the $\Ext$ functor, so it is $\BZ$-linear \cite[\S5,
n\degr3, prop.~6]{BourX}, \ie the $p$th power map in $K$ induces
multiplication by $p$ in $\Br(\K,K).$ Now we notice that the $p$th power
homomorphism induces an isomorphism of $F$ onto $\K,$ so it induces an
isomorphism of Brauer groups, and its composition with the embedding
$\K\subset F$ is the $p$th power map (and can be continued as the $p$th
power map to finite Galois extensions of $\K$), so the homomorphism
$\Br(\K)\to\Br(F)$ becomes multiplication by $p$ after composing with an
isomorphism, and it takes $p$-torsion to zero.

Thus we have constructed the required skew-field. Now we describe an 
example of a maximally central algebra over a perfect field that does not 
satisfy the condition \ref{:lincomb}) of Proposition \ref{bigcenter}. We take 
the field $\K$ as above, set $B=\K[\T t]/(\T t^{p^2}-x^p)$ and
$A=B\langle\sigma,u\rangle/(\sigma^p-(\T t-1), u^p-u-\T t,\sigma
u-(u+1)\sigma)$ (here $\sigma,u$ commute with $B$). This algebra is 
a cyclic crossed product in the sense of \cite[Section~6]{Az}. We see also 
that $a=\T t^p-x$ is a nilpotent of degree $p$ in $B$ and $B/(a)=F,$ so 
that $A/(a)$ is the skew-field we constructed. Let us note also that
 $B\subset Z(A).$ If we set $\sigma>u$ and take the degree-lexicographic 
order, then the initial terms of the defining relations of $A$ are
 $\sigma^p,$ $u^p,$ $\sigma u,$ and the monomials not divisible by them 
are only $u^i\sigma^j, 0\leq i,j <p,$ so $A$ is generated by $p^2$
monomials as a $B$-module. Hence the condition \ref{:maxc}) of
Proposition~\ref{bigcenter} is satisfied, and from this (from the fact that
the inequality cannot be strict) it follows that $Z(A)=B.$ If $A$ satisfied
\ref{:lincomb}), then, as the quotient of $A$ modulo the nilradical is a
skew-field, the sum would consist of just one summand, \ie
$A=S\otimes_{Z(S)}K$ for a simple $S$ and a commutative local $K$. Then
$A/\rad A=S\otimes_{Z(S)}K/\rad K,$ \ie the skew-field we constructed is
obtained from $S$ by extension of scalars from $Z(S)$ to $K/\rad K.$ Now,
the center of our skew-field is $F,$ and $[F:\K]=p,$ so $\K\subset
Z(S)\subset F$ implies that one of these inclusions is an equality. But the
skew-field we constructed cannot be obtained by extension of scalars from a
skew-field over $\K,$ thus $Z(S)=F.$ Therefore $F$ should be embedded into
$Z(A)=B$ as a $\K$-subalgebra with $1,$ but it cannot be embedded this way:
every preimage in $B$ of $t\in F$ has the form $\T t+ab$ for some $b\in B,$
and $(\T t+ab)^p=\T t^p+a^pb^p=\T t^p=x+a\ne x.$

\chapter{Proof of Theorem~\ref{stat1}}\label{sect1}

First we describe $M\mapsto\Fl M$ as a functor. Let $A$ be a
finite-dimensional associative algebra with 1, $M$ a finite-dimensional
 $A$-module. Let also $\{f_i\}$ be a $\K$-basis of $A,$ $A^l$ 
the $l$-fold direct product of $A$ with itself, considered as an algebraic
variety (an affine space),
$$
R=\K[A^l]=S((A^*)^l)=\K[\{x_{ij}\mid i\in[1,\dim_{\K}A],j\in[1,l]\}]
$$ 
its coordinate ring, where $S$ stands for the symmetric algebra of a
vector space and  $x_{ij}$ is the coefficient before $f_i$ in the $j$th
term of $A^l.$ We set 
$$
\T A = R\otimes_{\K}A,\quad
\Id_j=\sum_{i=1}^{\dim A} x_{ij}\otimes f_i\in \T A
$$
generic elements of $A,$ \ie the images of the identity operator
$\Id\in A^*\otimes_{\K}A$ under the mapping of $A^*\otimes_{\K}A$ into $\T
A,$ induced by the embedding of $A^*$ onto the $j$th summand of
$(A^*)^l\subset S((A^*)^l)=R,$ $\T M=R\otimes_{\K}M$ an $\T A$-module. 
(According to \cite[I, \S4]{Kro}, the coordinate construction of a generic 
element of an algebra $A$ as an element in $A\otimes_\K Q(R)$, $Q(R)$ 
being the quotient field of the polynomial ring $R$ with $l=1$, dates back 
to Kronecker.) Then 
\begin{gather*}
\Fl M=\T M/(\Id_1,\dots,\Id_l)\T M =\Fl A\otimes_{\T A}\T M=\\
=\Fl A\otimes_{\T A}(R\otimes_\K A)\otimes_A M=\Fl A\otimes_A M,
\end{gather*}
where $\Fl A=\T A/(\Id_1,\dots,\Id_l)\T A$ has the structure of a right 
$\T A$- (in particular $A$-) module induced from $\T A.$ 
This shows that $\Fl{\cdot}$ is an additive $\K$-linear right-exact functor
from the category of $A$-modules into the category of graded $R$-modules,
where the grading on $\Fl M$ is defined as follows. If we consider $M$ as a
vector space concentrated in degree 0, then the grading of $R$ and this 
grading of $M$ define a grading on $\T M,$ \wrt which $(\Id_1,\dots,\Id_l)\T M$ 
is a homogeneous $R$-submodule, so this grading induces a grading on the
quotient $\Fl M.$

We remark that the fiber over $(a_1,\dots,a_l)$ of the sheaf on $A^l$ 
corresponding to $\Fl M$ is the vector space $M/(a_1,\dots,a_l)M.$

\section{The $l=1$ Case}
We write simply $\Id$ for $\Id_1.$
\begin{lemma}\label{l1}
For every $A$-module $M$ the sequence of $R$-modules
$$0\to\T M\Idarrow \T M \to \Ff M\to 0$$ is exact.
\end{lemma}
\Proof
Note that $\Ff M = \T M/\Id\T M$ and that one has to check only the
exactness in the leftmost term, \ie the vanishing of the kernel of
$\Id$ on $\T M.$ 
We prove this by contradiction:
let $x\in\T M$ be a nonzero element in $\Ker\Id.$ Then $(\det\Id)x=0,$
and as $\T M$ is a free $R$-module and $R$ is a polynomial ring, one has
$\det\Id=0.$ But if we give the variables $x_i$ such values
$c_i\in\K$ that $\sum c_if_i=1_A,$ we shall obtain that
$\det\Id(c_1,\dots,c_{\dim A})=\det 1_M=1,$ hence $\det\Id\ne0$ ---
a contradiction.\qed

\vskip1em
\rem Actually we have just checked a particular case of the exactness
criterion of Prop.~\ref{prop.WMCE} and also shown that $\det\Id$ 
is a regular element in $\ann \Ff M.$
Therefore $\dim \Ff M\leqslant\dim R -1,$ but our exact sequence is the
minimal resolution of $\Ff M$ over $R,$ hence Prop.~\ref{dpft2} gives that 
$\Ff M$ is a Cohen-Macaulay module and $\pd_RM=1,$ so part~\ref{good}) of 
the theorem is proved.

\vskip1em
Next we prove the exactness of $\Ff{\cdot}.$ Set
$\CC=(0\to \T A \Idarrow \T A\to 0)$ to be a free resolution of 
$\Ff A$ over $\T A$ and $A$ (according to Lemma~\ref{l1}), then for any 
$A$-module $M$ $H_i(\CC\otimes_A M)=
\Tor_i^A(\Ff A, M),$ but
$$
\CC\otimes_A M =\left(0\to \T M \Idarrow \T M\to 0\right),
$$ 
whence it follows by Lemma~\ref{l1}
that for every $M$
$\Tor_1^{A}(\Ff A, M)=0.$
Now for every exact sequence of $A$-modules $0\to M\to N\to P\to 0$ one has
an exact sequence
$$
0=\Tor_1^A(\Ff A,P)\to\Ff A\otimes_A M\to\Ff A\otimes_A N\to
\Ff A\otimes_A P\to 0,
$$
as we needed.

\vskip1em
Now we prove that $\Ff{\cdot}$ is a fully faithful functor into the category 
of graded modules.
The minimal resolution of $\Ff M$ over $R$ is of the form
$$0\to\T M \Idarrow\T M \to \Ff M\to 0,$$
and any homogeneous homomorphism of
$R$-modules $\varphi\co\Ff M \to\Ff N$ can be extended to a graded
homomorphism of their minimal resolutions:
$$
\xymatrix{
 0\ar[r] &\T M \ar[r]^{\Id\cdot}\ar[d]_{\T g} 
&\T M \ar[r]\ar[d]_{\T f} &\Ff M\ar[r]\ar[d]_{\varphi} & 0 \\
0\ar[r] & \T N \ar[r]^{\Id\cdot} &\T N \ar[r] &\Ff N\ar[r] & 0.
}
$$

The homomorphisms $\T f$ and $\T g$ between free modules are given by 
matrices over $R,$ and as these homomorphisms are homogeneous of degree $0,$ 
the matrices are actually over $\K,$ \ie $\T f$ and $\T g$ result from
homomorphisms of vector spaces $f,g\co M\to N.$
The commutatrivity of the diagram yields $\T f\circ\Id=\Id\circ\T g.$
Substituting into $\Id$ the values $a_i,$ for which $\sum a_if_i=1_A,$
we obtain $f=g,$ and substituting all other values we obtain that
$f\in\Hom_A(M,N).$ As $\Ff f=\varphi,$ the map between the spaces of
homomorphisms is onto.
Since after choosing a basis in $\T M$ and projecting it to a 
$\K$-basis in $\Ff M/\mm\Ff M$ with $\mm=(x_1,\dots,x_{\dim A}),$
$\T f$ and $(R/\mm)\otimes\varphi$ are given by the same matrix, the matrix
of $f$ can be recovered from $\varphi$ and the map between the spaces of
homomorphisms is injective.

\rem For non-isomorphic $M$ and $N$ the modules $\Ff M$ and $\Ff N$ are not
isomorphic even as ungraded modules, as the homogeneous component of degree
0 of an isomorphism is itself an isomorphism.

\vskip1em
\rem\label{fullfaith} If we represent $\Bbbk[A]$ as the quotient of the 
ring $\Bbbk[A^l]$ modulo the ideal generated by the variables $x_{ij}$ for 
$j\geq 2,$ we see that
$$
F_l(M)\otimes_{\Bbbk[A^l]}\Bbbk[A]=F_1(M),
$$
since this multiplication does not change $\Id_1$, taking all other 
$\Id_j$ to zero. As the modules $F_l(M)$ are generated in degree $0$, an 
element in $\Homgr(F_l(M),F_l(N))$ is given by a matrix in 
$\Hom_\Bbbk(M,N)$, \ie the functors $F_l(\cdot)$ and the tensor product 
are injective on morphisms. As their composition is bijective on 
morphisms, we obtain the fact that \emph{$F_l(\cdot)$ is fully faithful as 
a functor into the category of graded modules for any associatvie algebra 
$A$ with $1$}.

\vskip1em
We prove part~\ref{three}). First we remark that the minimal number of
generators of the annihilator of a module is preserved under extension of
the base field, so until the end of the proof of part~\ref{three}) 
we assume $\K$ to be algebraically closed.

\begin{lemma}\label{l2}
If $M$ is a simple $A$-module, then $\ann \Ff M$ is a principal prime ideal.
\end{lemma}
\Proof As $M$ is a simple $A$-module, it is a module over a simple
summand of the semisimple quotient algebra of $A,$ \ie the module 
$\K^n$ over $M_n(\K),$ according to the theory of finite-dimensional
associative algebras (\cite[Kap.~13--14]{VDW}, \cite{CR}). Then
$\Ff M=R^n/\langle A\rangle$
with $A=(a_{ij}),$ $a_{ij}$ being independent variables. Now for such a
module one knows that,  $\det A$ being irreducible, the ideal $(\det A)$ is
prime and equals $\ann\Ff M$
($(\det A)\subset \ann\Ff M,$ $(\det A)$ is a prime ideal of height 1, and
$\pd_R \Ff M=1 \hence \hht\ann \Ff M \leqslant 1$). \qed

Let us prove part~\ref{three}) by induction on the length of the $A$-module $M.$ 
For $l(M)=1$ this is Lemma~\ref{l2}. If $l(M)>1,$ then there is a simple
submodule $N\subset M.$
Then part~\ref{good}) gives us an exact sequence
 $0\to\Ff N\to\Ff M\to\Ff P\to 0.$
As $l(P)<l(M),$ by the induction assumption $\ann\Ff P$ is a principal
ideal, say $(a).$ Take an $x\in\ann\Ff M,$ then $x$ annihilates $\Ff P.$
Hence $x=ay,$ where $y$ annihilates $a\Ff M.$ Now $a\Ff M$ is a 
submodule of $\Ff N;$ if $a\Ff M=0,$ then $x$ annihilates $\Ff M$ iff
 $x\in (a),$ so $\ann\Ff M=(a);$
but if $a\Ff M\ne0,$ then $\Ass a\Ff M\subset \Ass
\Ff N.$ But $\Ass \Ff N=\{(p)\}$ for $(p)=\ann\Ff N,$ because 
by part~\ref{good}) $\Ff N$ is Cohen-Macaulay, so $\Ass\Ff N$ coincides 
with the set of
the minimal primes over $\ann\Ff N,$ and $\ann\Ff N$ is prime. Then, as all
the assosiated primes of  $a\Ff M$ contain its annihilator, one has
$\ann a\Ff M\subset(p)=\ann\Ff N,$ but $a\Ff M$ is a submodule of $\Ff N,$
therefore $\ann\Ff N\subset\ann a\Ff M.$ Hence $\ann a\Ff M = (p).$
Thus $\ann\Ff M=a\ann a\Ff M=(ap).$ \qed

\vskip1em
\rem Let us look at this construction under more general settings. Assume
that we consider representations of a vector space 
$V=\langle f_1,\dots,f_d\rangle_{\K}$ with $1=f_1,$ \ie with an element
that acts as identity in any representation. Then we take for $A$ a free 
algebra $\K\langle f_2,\dots,f_d\rangle$ and consider the representations of 
$V$ as $A$-modules. In this case the description of the functor $\Ff{\cdot}$
remains the same %Led Zeppelin rulez!
$(R=\K[V]=\K[x_1,\dots,x_d])$ and the same proofs are valid for two
first parts of Theorem~\ref{stat1}, only the fact that $\CC$ is a 
resolution for $\Ff A$ does not follow from Lemma~\ref{l1}, but is obtained by 
considering
the lowest homogeneous component of $\Id v$ \wrt the grading of the free
algebra $A,$ as this component equals the lowest component of $v$ times $x_1.$

The presence of the identity is essential for the proofs of parts 1) and
2): for
example, if one takes the adjoint representation of a non-Abelian
finite-dimensional Lie algebra and constructs the corresponding module, it
\emph{won't} be Cohen-Macaulay,
 but will be of Krull dimension $\dim R$ and of projective
dimension $\geqslant 2$: let $L$ be the Lie algebra, then
$\Id\in\T L,\ \Ff L=\T L/[\Id,\T L],$ hence
$$
\T L\myto{[\Id,\cdot]}\T L\to\Ff L\to0
$$
is the beginning of the minimal resolution of $\Ff L$ over $R$ and the
resolution continues to the left, because $\Id\in\Ker[\Id,\cdot].$
But $\dim \Ff L = \dim R,$ for, as observed by the referree of \cite{art1},
$\mathop{\rm supp} \Ff L$ in $\K^{\dim L}\cong L$ equals
$$
\{x\in L \mid L/[x,L]\ne 0\}=L,
$$ 
because
$\forall x\in L$ $\dim_{\K} L - \dim_{\K} [x,L] = \dim_{\K} \Ker [x,
\cdot]>0.$

\section{The $l>1$ Case}\label{lg1}
First of all, we remark that a maximally central algebra is by definition 
a direct sum of equidimensional maximally central ones, the modules over 
it are direct sums of modules over the summands and the same is true for short 
exact sequences and homomorphisms (\ie for categories). In particular, 
indecomposable modules are modules over one of the summands. So it 
suffices to prove the theorem for equidimensional maximally central
algebras: it suffices as well to check the exactness of the functor over
each of the summands.

Then come some general remarks on the extension of scalars.
\begin{lemma}\label{ACF}
\begin{enumerate}
\item 
The functor $\Fl{\cdot}$ commutes with the extension of scalars.
\item\label{ACF-CM}
The module $\Fl M$ is Cohen-Macaulay after the extension of scalars iff 
it was Cohen-Macaulay before the extension.
\item
Statement~\ref{two}) of Theorem~\ref{stat1} is satisfied after the extension of
scalars iff it was satisfied before the extension also.
\end{enumerate}
\end{lemma}
\Proof
\begin{enumerate}
\item 
Evident from the construction.
\item 
We note that for every graded $R$-module $D$ extension of scalars 
preserves the minimal graded resolution of $D$ over $R,$ hence the 
projective dimension and the Krull dimension (as the Hilbert function of a
module is the Euler characteristic of the resolution \wrt the Hilbert
function) of $D,$ so the presence or absence of Cohen-Macaulayness as well
(according to prop.~\ref{dpft2} and \ref{grad.CM}, $D$ is Cohen-Macaulay 
$\Leftrightarrow \pd D+\dim D=\dim R$).
\item 
$\Hom$ and homology commute with the extension of scalars.\qed
\end{enumerate}

Thus it suffices to prove the theorem for an equidimensional maximally 
central algebra and after passing to the algebraic closure of the base 
field. So by Remark~\ref{barequidim} it suffices to prove the theorem for
algebras of the form $M_n(K),$ $K$ being commutative, over an algebraically
closed field $\K,$ and this will occupy the rest of the section.

Note that for $A=M_n(K)$ the tensor multiplication by the $(A,K)$-
and $(K,A)$-bimodule $K^n$ defines an equivalence between the categories of
$A$-modules and $K$-modules (Morita-equivalence). We now describe the
functor which associates to an $A$-module $M$ the minimal resolution of
$\Fl M$ over $R.$

\refstepcounter{lemma}{\bf Notations \thelemma.}\label{Not}
The ring $R$ can be described as 
$\K[\{x_{\alpha ij} \mid \alpha\in[1,\dim_{\K} K],\allowbreak{}i\in[1,n],
\ j\in[1,ln]\}].$ We set $R[T]=R[\{T_{ij} \mid i\in[1,n],\ j\in[1,ln]\}],$
$\varphi=(T_{ij})$ an $(n\times ln)$-matrix,
$C$ the number 1 Eagon-Northcott complex corresponding to $\varphi$
(Def.~\ref{EN}),
$\T K=K\otimes_{\K}R,\ \Id^K_{ij}=\sum_{\alpha}f_{\alpha}x_{\alpha ij}
\in \T K,$ $f_{\alpha}$ being a $\K$-basis for $K.$
Let $\KK$ denote the $(\T K,\T A)$-bimodule $\T K^n.$
The operation of $T_{ij}$ on $\KK$ by means of the commuting ($K$ being
commutative) endomorphisms $\Id^K_{ij}$ defines an $R[T]$-module structure
on $\KK.$ Set $\CC=C\otimes_{R[T]}\KK,$ a complex of
projective (as direct sums of $\KK$) right $\T A$-modules.

Then $\Fl M=H_0(\CC\otimes_{A} M)$ for every $M$: is is 
$\Coker(\T\M^{ln}\to\T\M^n)$ with $\T\M=\KK\otimes_{\T A}\T M,$
and \wrt $R$-bases this homomorphism is given by a block matrix
$\varphi=(\Id^K_{ij}).$ The complex $\CC\otimes_{A} M=
C\otimes_{R[T]}\T\M$ is a minimal complex of length $(l-1)n+1$
and it is the minimal resolution of $\Fl M$ and $\CC$ is a projective
resolution of $\Fl A$ as a right $\T A$- (and $A$-) module by virtue of 
the following lemma:

\begin{lemma}\label{Mn}
For any $M$ and any $i>0$ $H_i(\CC\otimes_A M)=0.$
\end{lemma}
\Proof According to the exactness criterion for the Eagon-Northcott complex
(Subsection~\ref{WMCE}), it suffices to prove that the ideal $I_n(\varphi)$ 
of the maximal minors of $\varphi$ contains an
$\T\M$-sequence of legth $(l-1)n+1.$ This sequence, as we show, consists of
the minors corresponding to $n$ columns in succession. $\T\M$ is a
Cohen-Macaulay 
$R$-module, hence a Cohen-Macaulay $R[T]$-module (as a graded module over 
a graded ring,
where the new variables $T_{ij}$ also have degree 1), because its Krull
dimension, determined from the Hilbert function of the grading, is the same
over both rings, and a maximal homogeneous $\T\M$-sequence in $R$ is also a
maximal homogeneous $\T\M$-sequence in $R[T].$ Therefore (Prop.~\ref{regCM}),
a regular $\T\M$-sequence is a sequence, after factoring which out the
dimension of the module decreases by the length of this sequence.

 $\T\M$ is a free $R$-module, and the quotient module of $\T\M$ by our
sequence is the quotient by the columns of the matrices which are the sums
of products of $\Id^K_{ij}$ corresponding to the minors. Let $f_1$ be the
identity of $K.$ We order the variables $x_{\alpha ij}$ the following way:
$$
\renewcommand{\land}{\mathbin{\&}}
x_{\alpha ij}>x_{\beta lm}\rightleftharpoons (i<l)\lor(i=l\land j<m)\lor
(i=l\land j=m\land \alpha<\beta)
$$ 
and take the corresponding degree-lexicographic order in $R$ and the 
corresponding ``term over position" order in $\T\M$ (for some order of the 
basis elements of $\T\M$).

What are the initial terms of our columns of relations \wrt this order?
Every monomial in every element of the matrix corresponding to
$T_{i_1j_1}\cdots T_{i_nj_n}$ is
$x_{\alpha_1i_1j_1}\cdots x_{\alpha_ni_nj_n}$ for some $\alpha.$
The greatest of these monomials is $x_{1i_1j_1}\cdots x_{1i_nj_n}.$
The fact that the substitution $x_{2ij}=\dots=x_{(\dim K)ij}=0$ results in
the matrix $x_{1i_1j_1}E\cdots x_{1i_nj_n}E,$ $E$ being the identity matrix,
shows that this monomial appears in all diagonal elements, therefore in all
columns. So the initial terms of the columns of the matrix corresponding
to $T_{i_1j_1}\cdots T_{i_nj_n}$ are the columns of the matrix
$x_{1i_1j_1}\cdots x_{1i_nj_n}E.$

Now we see, that if we add two products corresponding to different sets
$(i_1j_1,\dots,i_nj_n),$ then the initial terms are different and do not
cancel out, so the initial terms of the columns of the matrix corresponding
to the minor from the columns $j+1,\dots,j+n$ are the columns of the matrix
$$\max_{\sigma\in S_n} x_{11(j+\sigma(1))}\cdots x_{1n(j+\sigma(n))}E=
x_{11(j+1)}\cdots x_{1n(j+n)}E.$$

Indeed $x_{11(j+1)}$ is the greatest variable occuring in the initial
monomials, $x_{12(j+2)}$ is the greatest one occuring in the monomials
containing $x_{11(j+1)}$ etc.. Thus the initial terms of our relations are
the columns of the matrices
$y_1E,y_2E,$ $\dots,y_{(l-1)n+1}E$ with $y_i=x_{11i}\cdots x_{1n(i+n-1)}.$

As any two leading terms either have different basis vectors or depend 
from disjoint sets of variables, there are no critical pairs, so 
the columns we consider constitute a Gr\"obner basis of relations. Thus 
the Krull dimension of the quotient modulo the columns is the same as that 
of the quotient modulo the leading terms, \ie of the module 
$\T\M/(y_1,\dots,y_{(l-1)n+ 1})\T\M.$ Since $\T\M$ is a free $R$-module 
and $y_i$ depend on disjoint sets of variables in $R,$ they form an 
$\T\M$-sequence, so the Krull dimension of the quotient modulo the leading 
terms, and therefore of the quotient modulo the minors, underwent the 
required decrease. \qed

Thus our complex is really the minimal resolution of
$\Fl M,$ and also $\forall M\ \Tor_1^{A}(\Fl A, M)=0.$
Now one can show the exactness of $\Fl{\cdot}$ in the same fashion as for $l=1$:
\begin{gather*}
0\to M\to N\to P\to 0\text{ is exact }\hence \\
0=\Tor_1^A(\Fl A,P)\to\Fl A\otimes_A M\to\Fl A\otimes_A N\to 
\Fl A\otimes_A P\to 0
\end{gather*}
is exact. The full faithfullness of $\Fl\cdot$ follows from 
Remark~\ref{fullfaith}, so part~\ref{two}) is proved.\qed
\vskip1em
 Let us prove the Cohen-Macaulayness. For this one might similarly construct
a regular sequence in $\ann\Fl M,$ but we pursue another way: induction
on the length of $M.$ If $\K$ is algebraically closed and $M$ is simple, then 
$\Fl M$ is the quotient of a free module by the columns of the generic
$(n\times ln)$-matrix, and its Cohen-Macaulayness is well-known
(\cite[Appendix 2.6]{Ei}, \cite{BV}). Further, let $l(M)$ be greater than $1$
and $0\to M_1\to M\to M_2\to 0$ be an exact sequence of $A$-modules with
$l(M_i)<l(M).$ Then by the induction assumption $\dim M_i=
\depth M_i=\dim R -((l-1)n+1)$ and by the features of depth and of Krull
dimension
(Prop.~\ref{dpft}; \cite[chap.~III, B.1, chap.~I, C.1, prop.~10]{Se})
$\dim M \leqslant\max_i\dim M_i,\depth M\geqslant\min_i \depth M_i,$ that
is, $\dim M=\depth M$ and $M$ is also Cohen-Macaulay. \qed

\rem Part~\ref{three}) emerged from the hope that $R/\ann\Fl M$ is
Cohen-Macaulay and $\Fl M$ is a maximal Cohen-Macaulay module over it 
(\ie of maximal
dimension). But this is not the case already for $l=2$ and for the standard
representation of the algebra of diagonal $2\times2$-matrices, when this
ring has the form $\K[x_1,x_2,y_1,y_2]/(x_iy_j)$ and Krull dimension 2, and
after factoring out the regular element $x_1-y_1$ the element $x_1$ will be
annihilated by all the variables.

\chapter{Proof of Theorem~\ref{stat2}}

Throughout this chapter we assume $l>1.$

\section{The Case of Cohen-Macaulayness}

We prove the strongest claim at once, since it is easily reduced to
the case of an algebraically closed field.

\begin{satz}
If for some $l>1$ and some $A$-module $M$ $F_l(M)$ is Cohen-Macaulay, then
$A/\ann M$ is an equidimensional maximally central algebra.
\end{satz}
\Pf Note that we can assume the field to be algebraically closed,
as the Co\-hen-Ma\-cau\-lay\-ness of $F_l(M)$ does not depend on this 
(part~\ref{ACF-CM}) of Lemma~\ref{ACF}). So by Remark~\ref{barequidim} we 
have to show that $A/\ann M$ is a direct sum of matrix algebras of the 
same rank over commutative ones. Passing to the quotient of $A$ modulo 
$\ann M$ we can assume $\ann M=0$. So \emph{in the sequel of the section 
we assume that $\K$ is algebraically closed and $\ann M=0$ ($M$ is a 
faithful module)}.

\begin{lemma}\label{firstfact}
$A$ is a direct sum of algebras $A_i$ the semisimple quotients of which 
are simple, and these quotients are matrix algebras of the same rank.
\end{lemma}
\Pf
Choose an embedding of $A/\rad A=M_{n_1}(\K)\oplus\dots\oplus M_{n_k}(\K)$ 
into $A$ as a subalgebra with $1$ (you can choose one over an 
algebraically closed field by the Wedderburn--Malcev Theorem, see 
Prop.~\ref{WeddMal}), choose a composition series $0=M_0\subset
M_1\subset\dots\subset M_m=M$ in $M$ and choose a $\K$-basis
$e_1,\dots,e_r$ in $M$ that conforms to these choices, \ie that $M_i=\langle
e_1,\dots,e_{k_i}\rangle_\K$ and that $\langle 
e_{k_{i-1}+1},\dots,e_{k_i}\rangle_\K$ are simple submodules over the 
subalgebra $A/\rad A$. In the decomposition 
$A=M_{n_1}(\K)\oplus\dots\oplus M_{n_k}(\K)\oplus\rad A$ decompose $\rad 
A$ further into isotypic components as a bimodule over the semisimple part 
and choose a $\K$-basis in $A$ that conforms to the resulting 
decomposition (matrix elements being the basis of the semisimple part). 
Then the matrices $\Id_j$ become blockwise upper triangular, with the
diagonal occupied by square blocks of independent variables corresponding
to the simple quotients of $A$ over which the corresponding composition
factors are simple modules, and with linear forms in variables
corresponding to the radical of $A$ above the diagonal. One can say more:
the intersection of a row that has a simple quotient of type $P$ on the
diagonal and a column that has a simple quotient of type $Q$ contains
linear forms in the variables corresponding to the $(P,Q)$-isotypic
component of the radical.

We shall illuminate the behavior of the matrices $\Id_j$ under the 
factorizations we do in the proof by the example of a module of length $4$ 
with successive composition factors of types $P,$ $P,$ $Q,$ $P$:
$$
\Id_j=
\begin{pmatrix}
P & (P,P) & (P,Q)& (P,P)\\
& P& (P,Q)& (P,P)\\
&& Q &(Q,P)\\
&\raisebox{0pt}[0pt]{\makebox[0em][r]{\Huge $0$\ }} && P
\end{pmatrix}.
$$

Let $M_{n_1}(\K)$ be the simple quotient of the algebra to which 
corresponds the last composition factor $P$ of our composition series (\ie 
the simple quotient of our module). Take the quotient of $F_l(M)$ modulo 
the sequence of the variables corresponding to the radical of the algebra 
and then by the sequence of the variables corresponding to the remaining 
simple quotients of the algebra in elements $\Id_2,\dots,\Id_l.$ Then take 
the quotient of the resulting module modulo the sequence of the variables 
that correspond to the remaining simple quotients in $\Id_1$ and do not 
occupy its principal diagonal and modulo the sequence $y_s-1,$ $y_s$ running 
over the variables that correspond to the other simple quotients in 
$\Id_1$ and do occupy its principal diagonal. Then we obtain a module $X'$ over 
the ring $R'$ of polynomials in the remaining variables (\ie corresponding 
to the simple quotient chosen), for which the presentation matrix has 
only blocks of variables that correspond to this quotient on the diagonal, 
zeros above it, and the blocks that correspond to other simple quotients 
of $A$ are turned to identity matrices in $\Id_1$ and to zero matrices in
$\Id_2,\dots,\Id_l$:
$$
\Id_1'=
\begin{pmatrix}
P & 0 & 0& 0\\
& P& 0& 0\\
&& 1 & 0\\
&\raisebox{0pt}[0pt]{\makebox[0em][r]{\Huge $0$\ }} && P
\end{pmatrix},\quad
\Id_j'=
\begin{pmatrix}
P & 0 & 0& 0\\
& P& 0& 0\\
&& 0 &0\\
&\raisebox{0pt}[0pt]{\makebox[0em][r]{\Huge $0$\ }} && P
\end{pmatrix}.
$$
So we can cross out the rows and columns in the presentation matrix
containing the identity matrices and realize that we obtain a direct sum of
several copies of $F^{R'}_l(P)$ over the polynomial ring in the remaining
variables.

The Krull dimension of $F^{R'}_l(P)$ (the quotient modulo the columns of 
an $(n_1\times ln_1)$-matrix of independent variables) is known to be
$\dim R'-(l-1)n_1-1$ \cite{BV}. Since $F_l$ is right-exact (see the 
beginning of Chap.~\ref{sect1}), $F^R_l(P)$ is a quotient of $F_l(M).$ 
Thus $\dim F_l(M)\geq\dim F^R_l(P)=\dim R-(l-1)n_1-1.$ As $F_l(M)$ is 
Cohen-Macaulay, its dimension decreased with the factorization by at most 
the length of the sequence factored out, namely $\dim R-\dim R',$ and in 
case of equality this sequence is regular (Prop.~\ref{regCM}). Therefore
$\dim F_l(M)=\dim F^R_l(P)$ and the sequence is regular.

Now we consider the quotient of $F_l(M)$ modulo a part of this sequence: we
take only those variables corresponding to nilradical in
$\Id_2,\dots,\Id_l$ that correspond to isotypic components not isomorphic
to direct sums of $P$ as left modules over the semisimple part of $A.$
Then $\Id_1$ and diagonals of other matrices change as above, and above 
the diagonal of $\Id_2,\dots,\Id_l$ the rows corresponding to composition 
factors of type $P$ remain the same, while other rows vanish:
$$
\Id_1''=
\begin{pmatrix}
P & 0 & 0& 0\\
& P& 0& 0\\
&& 1 & 0\\
&\raisebox{0pt}[0pt]{\makebox[0em][r]{\Huge $0$\ }} && P
\end{pmatrix},\quad
\Id_j''=
\begin{pmatrix}
P & (P,P) & (P,Q)& (P,P)\\
& P&(P,Q)& (P,P)\\
&& 0 &0\\
&\raisebox{0pt}[0pt]{\makebox[0em][r]{\Huge $0$\ }} && P
\end{pmatrix}.
$$
If we induce a filtration on this quotient $X''$ by the composition series
of $M,$ the adjoint factors of this filtrarion are either quotients of
$F^{R''}_l(P)$ over the ring in the remaining variables, as the relations
contain at least the diagonal blocks, or zeroes, if the diagonal block is
an identity matrix. Moreover, the last factor is exactly $F^{R''}_l(P).$ So
the Krull dimension of $X''$ equals $\dim R''-(l-1)n_1-1,$ and the same
argument as in the previous paragraph shows that the sequence is regular
and thus (Prop.~\ref{regCM}) $X''$ is Cohen-Macaulay.

As in the previous factorization, we can remove rows and columns 
occupied by identity matrices in $\Id_1$ and obtain a presentation of 
$X''$ by a matrix of linear forms. Thus $X''$ is a graded module and the 
passage from it to $X'$ is factoring out a homogeneous regular sequence of 
degree $1.$ Hence the Hilbert series of $X'$ can be obtained from the one 
of $X''$ by multiplying by $1-t$ raised to the power equal to the length of 
the sequence (see Section~\ref{inv}). So $X''$ and $X'\otimes_{R'}R''$ 
have the same Hilbert function. If we consider the filtration on 
these modules induced by the composition series of $M$, then, as remarked 
in the previous paragaraph, the adjoint factors for $X''$ are quotients of 
those for $X'\otimes_{R'}R'',$ so there is actually no further 
factorization. Suppose that after a factor of type $P$ we have a factor of 
another type $Q$ in the composition series (as in the example).
Then the part of the presentation matrix obtained from $\Id_2$ has some 
forms in the variables that correspond to the $(P,Q)$-isotypic component 
of the nilradical to the right of the block corresponding to the first 
factor and above the block corresponding to the second factor, with zeroes 
in the place of the second block on the diagonal and below:
$$
\Id_2''=
\begin{pmatrix}
P & (P,P) & (P,Q)& (P,P)\\
& P&(P,Q)& (P,P)\\
&& 0 &0\\
&\raisebox{0pt}[0pt]{\makebox[0em][r]{\Huge $0$\ }} && P
\end{pmatrix}.
$$
Thus if these forms are nonzero, they give an additional factorization of
the adjoint factor, which cannot be. Therefore they are zero and if we
transpose the corresponding groups of the basis vectors of $M$, we can put
the composition factor of $M$ isomorphic to $P$ after the one isomorphic to
$Q$:
$$
\Id_j=
\begin{pmatrix}
P & (P,Q) & (P,P)& (P,P)\\
& Q& 0& (Q,P)\\
&& P &(P,P)\\
&\raisebox{0pt}[0pt]{\makebox[0em][r]{\Huge $0$\ }} && P
\end{pmatrix}.
$$

Thus we can suppose that in the composition series of $M$ all the factors 
isomorphic to $P$ go after other factors and form a quotient module $M_P.$
Then if we factor out a sequence that does the same with the variables 
corresponding to the semisimple part and all the variables corresponding 
to the nilradical but for those corresponding to the $(P,P)$-isotypic 
component, $F_l(M)$ becomes $F_l^{R'''}(M_P)$ over the polynomial ring in 
the remaining variables:
$$
\Id_1'''=
\begin{pmatrix}
1 & 0 & 0& 0\\
& P & (P,P)& (P,P)\\
&& P & (P,P)\\
&\raisebox{0pt}[0pt]{\makebox[0em][r]{\Huge $0$\ }} && P
\end{pmatrix},\quad
\Id_j'''=
\begin{pmatrix}
0 & 0 & 0& 0\\
&P& (P,P)& (P,P)\\
&& P &(P,P)\\
&\raisebox{0pt}[0pt]{\makebox[0em][r]{\Huge $0$\ }} && P
\end{pmatrix}.
$$

The same filtration shows that
$$
\dim F_l^{R'''}(M_P)=\dim F_l^{R'''}(P)=\dim R'''-(l-1)n_1-1,
$$ 
so we were factoring out a regular sequence and $F_l^R(M_P)=
F_l^{R'''}(M_P)\otimes_{R'''}R$ is a Cohen-Macaulay module of the same 
dimension as $F_l(M).$ Therefore, according to the behavior of the depth 
(Prop.~\ref{dpft}) and Krull dimension \cite[chap.~III, B.1, chap.~I,
C.1, prop.~10]{Se} in short exact sequences, the leftmost module 
in the exact sequence
$$
0\to Y\to F_l(M)\to F_l(M_P)\to 0
$$
is Cohen-Macaulay of the same dimension. We show that this sequence is 
split. Because of Yoneda's interpretation of $\Ext$ \cite[\S7
n\degr3]{BourX} it will suffice for this to show that $\Ext^1_R(F_l(M_P),Y)=0.$
It is known (Prop.~\ref{dpft}) that
$$
\min\{i|\Ext^i_R(M,N)\ne 0\}=\depth(\ann M, N).
$$
It is enough to show that in our case this depth is at least 2. As the 
modules under consideration are Cohen-Macaulay and the ring is a 
polynomial ring (so that for a (prime, and hence any) ideal $I$ one has 
$\hht I+\dim R/I=\dim R$ \cite[chap.~III, D.3, prop.~15]{Se}, and the same
is true for any of its quotient rings), the depth equals $\hht(\ann Y+\ann
F_l(M_P))-\hht \ann Y$ \cite[Exercise~18.4]{Ei}, and all the annihilators
can be replaced by their radicals. As our modules have the same dimension,
$\hht\ann Y=\hht\ann F_l(M_P).$ Thus we have to show that
$$
\hht(\ann Y+\ann F_l(M_P))\geq\hht(\ann F_l(M_P))+2.
$$ 
The filtrarion induced by the composition series shows that the radical of 
the annihilator of $F_l(M_P)$ equals the radical of the annihilator of 
$F_l(P)$, namely, the ideal of the maximal minors of the matrix obtained 
by writing the blocks corresponding to $P$ in $\Id_j$ one after another.
It also shows that the radical of the annihilator of $Y$ contains the 
product of the ideals of maximal minors of matrices obtained similarly for 
other simple quotients of $A$. Such an ideal for an $n_i\times ln_i$ matrix 
is of height $(l-1)n_i+1\geq (2-1)1+1=2,$ thus their product is also of 
height $\geq 2$ and is generated by polynomials in variables not 
involved in the generators of the radical of $\ann F_l(M),$ so, if we add 
this product to the radical of $\ann F_l(M)$, its height increases by at 
least $2,$ as required.

Remark~\ref{fullfaith} says that $F_l$, though not always exact, is always 
fully faithfull as a functor into the graded category. Thus if $F_l(M)\to
F_l(M_P)$ is a split epimorphism (as the homogeneous component of degree 
$0$ of a left inverse to the projection is itself a left inverse, it 
doesn't matter whether the epimorphism is split in the graded or in the 
usual category), $M\to M_P$ is also a split epimorphism. Hence 
$M=M_P\oplus M'$ with $M'$ having no composition factors of type $P.$ Thus 
we can choose a composition series in $M$ in which the last quotient has 
some other type $Q$. Repeating the previous argument for this composition 
series, we see that $M=M_P\oplus M_Q\oplus M''$ with $M''$ containing no
composition factors of types $P$ and $Q.$ Induction on the number of
composition factors for which the corresponding ``isotypic components'' are
direct summands of $M$ yields us the conclusion that $M=\bigoplus_P M_P$ is
a direct sum of modules having only one type of composition factors each.
The matrices of the representaion of $A$ in a $\K$-basis of $M$ conforming
to this decomposition and further decompositions as in the beginning of the
proof are blockwise diagonal, and, as this representation is faithful, the
idempotents that correspond to the identities of the summands $M_{n_i}(\K)$
and are represented by matrices having one block identity and others zero lie
in the center of $A$ and decompose it into a direct sum of algebras
$\rho_{M_P}(A),$ each having only one simple module.

In the beginning of the proof we saw that $\dim F_l(M)=\dim F_l(P).$ Now 
we see that it is true for every $P,$ and as the latter dimension equals 
$\dim R-(l-1)n_i-1$ and we have $l>1,$ all the $n_i$ are equal. \qed

Now we can apply Lemma~\ref{matr} and write $A=M_n(K)$ for $K$ a direct 
sum of algebras with semisimple quotients equal to $\K$ (because the 
central idempotents that describe the decomposition of $A$ lie in $K$).
In Notations~\ref{Not} Morita-equivalence of the categories of modules over
$A$ and $K$ associates to a faithful $A$-module $M$ a faithful $K$-module 
$\M.$  Let $l'=(l-1)n+1\geq (2-1)1+1=2.$ A composition series of $\M$ 
shows that $\dim F_{l'}(\M)=\dim \K[K^{l'}]-l'.$ 

\begin{lemma}
$F_{l'}(\M)$ is a Cohen-Macaulay module, $\M$ being regarded as a $K$-module.
\end{lemma}
\Pf
In the notations~\ref{Not} $\Fl M$ is the quotient modulo the columns of an 
$n\times ln$ block matrix $(\Id^K_{ij}).$ Let $f_1$ be the identity of $K.$ 
Consider the sequence
$$
\mathbf{x}=(\{x_{1ij}-\delta_{ij}\}_{i,j=2}^n,\linebreak[0]
\{x_{\alpha ij}\}_{i,j=2,}^n{}_{\alpha=2}^{\dim_{\K}K},\linebreak[0]
\{x_{\alpha ij}\}_{i\geq2,j\not\in[2,n]},
\linebreak[0] \{x_{\alpha 1j}\}_{j=2}^n),
$$ 
\ie the one modulo which the matrix above becomes
$$
\left(
\begin{array}{ccccccc}
\Id^K_{11}&\multicolumn{3}{c}{\mbox{\Large 0}}&\Id^K_{1,n+1}&\dots&
\Id^K_{1,ln}\\
0&1&\multicolumn{2}{c}{\mbox{\Large 0}}\\
\vdots&&\ddots&&\multicolumn{3}{c}{{\mbox{\Huge 0}}}\\
0&\multicolumn{2}{c}{\mbox{\Large 0}}&1
\end{array}
\right).
$$

The quotient modulo the columns of the resulting matrix is easily seen to 
be $\Fll{\M}$, and one sees also that $\K[(M_n(K))^l]/(\mathbf x)=\K[K^{l'}]$, 
\ie $\Fll{\M}=\Fl M/\mathbf x\Fl M.$ We need to prove that $\Fll{\M}$ is 
Cohen-Macaulay, and since $\Fl M$ is Cohen-Macaulay by hypothesis, it 
suffices to prove by Prop.~\ref{regCM} that $\dim \Fll{\M}=\dim\Fl M -
l(\mathbf x).$
Now, $\mathbf x$ is a regular sequence of length $\dim \K[(M_n(K))^l]-
\dim \K[K^{l'}]$ in $\K[(M_n(K))^l]$, and we have remarked above that the 
dimensions $\dim\Fl
M=\dim\K[(M_n(K))^l]-l'$ and $\dim\Fll{\M}=\dim\K[K^{l'}]-l'$ are as 
required, therefore $F_{l'}(\M)$ is Cohen-Macaulay. \qed

\begin{lemma}
Let $K$ and $\M$ be as above. Choose an embedding of $K/\rad K=\K^k$ into 
$K$ as a subalgebra containing $1_K,$ choose a composition series 
$0=\M_0\subset\M_1\subset\dots\subset \M_m=\M$ in $\M$ and choose a
$\K$-basis $e_1,\dots,e_m$ for $\M$ conforming to these choices, \ie so 
that $\M_i=\langle e_1,\dots,e_i\rangle_\K$ and $\K e_i$ are simple 
submodules over the subalgebra $K/\rad K$. Choose a $\K$-basis for $K$ 
conforming to the decomposition $K=\K^k\oplus\rad K$ and a monomial order 
in $R=\K[K^{l'}].$ Let $e_1<e_2<\dots<e_m$ and introduce the ``position 
over term'' monomial order in $R^m$. Then the columns of the matrix
$\Id^K_1|\dots|\Id^K_{l'}$ form a Gr\"obner basis of the submodule they 
generate.
\end{lemma}
\Pf
Note that under these choices $\Id^K_i$ is upper triangular, the diagonal 
occupied by the variables corresponding to simple quotients of $K$ and the 
space above it by linear forms in the variables corresponding to the 
radical of $K,$ so that the initial terms of the columns of this matrix 
are obtained by substituting zero into all the variables corresponding to 
the radical.

In analogy with the proof of Lemma~\ref{firstfact}, consider the quotient of 
$F_{l'}(\M)$ modulo a sequence of the variables corresponding to the basis 
of $\rad K.$ It is a regular sequence in $\K[K^{l'}]$ and $\K[(K/\rad 
K)^{l'}]$ is the quotient of the ring by it, while the quotient of 
$F^K_{l'}(\M)$ modulo this sequence is $F^{K/\rad K}_{l'}(\M)$ for the 
restriction of $\M$ to the semisimple quotient embedded as a subalgebra.
The standard argument involving a composition series and the 
right-exactness of $F_{l'}$ shows that the Krull dimensions of both 
modules equal the dimensions of the corresponding rings less $l'$. Thus 
this sequence is $F^K_{l'}(\M)$-regular, so the Hilbert series of the 
quotient is obtained from the Hilbert series of $F^K_{l'}(\M)$ through
multiplication by $1-t$ raised to the power equal to the length of the 
sequence. Therefore if we tensor $F^{K/\rad K}_{l'}(\M)$ over $\K$ with 
the polynomial ring in the variables corresponding to the radical, we
obtain a module with the same Hilbert function as $F^K_{l'}(\M).$ But the
module we obtain is the quotient of $R^m$ modulo the initial terms of the
relation columns for $F^K_{l'}(M),$ whence these columns are a Gr\"obner
basis. \qed

\begin{lemma}
The algebra $K$ considered above is commutative.
\end{lemma}
\Pf
Consider the commutator of two generic elements. Its columns lie in the 
submodule of relations and should be reducible to zero, but they depend 
only on the variables corresponding to the nilradical, as the matrices 
multiplied by the variables corresponding to $\K^n$ are central 
idempotents, whereas the initial terms of the Gr\"obner basis are 
divisible by a variable corresponding to the semisimple quotient. So the 
commutator equals zero.
\qed

We have proved the proposition. Now we prove the remainder of the theorem. 
If $F_l$ takes all $A$-modules to Cohen-Macaulay ones, then, applying the 
proposition to the left regular representation of the algebra (which is 
faithful) we see that the algebra is equidimensional maximally central. It 
also follows from the proposition that if $F_l(M)$ is Cohen-Macaulay for 
an indecomposable $A$-module $M$, then all the composition factors of $M$ 
are isomorphic. Then it follows from the description of the decomposition of 
finite-dimensional algebras into a direct sum in Section~\ref{decomp} that 
$A$ is a direct sum of algebras with simple semisimple quotients. So for 
every summand all the indecomposable modules are taken into Cohen-Macaulay 
modules of the same dimension (as for the simple module), so all modules 
are taken into Cohen-Macaulay ones, so the summands, and therefore the 
whole $A$, are maximally central.

\section{The Case of an Exact Functor}

First we prove that the functor $F_l(\cdot)$ remains exact over $\B\K$, 
the algeraic closure of $\K$. We have seen at the beginning of 
Chapter~\ref{sect1} that $F_l(\cdot)$ is the tensor multiplication over 
$A$ by the right $A$-module $F_l(A)$. But the flatness of a module is known 
to be equivalent to preserving the exactness of sequences of finitely 
generated modules under tensoring with this module \cite[\S4, n\degr6, 
th\'eor\`eme 2]{BourX}. So the exactness of $F_l(\cdot)$ is equivalent to
 $F_l(A)$ being a flat right $A$-module. Now, this condition is preserved 
under extension of scalars and under passing to induced modules in 
general, which is immediate from the associativity of the tensor product.

So we can assume $\K$ to be algebraically closed. Now we show that 
non-isomorphic simple $A$-modules form no nontrivial extensions.

\begin{lemma}
Let
\begin{equation}
\label{eq1}
0\to M\to N\to P\to 0
\end{equation}
be a non-split exact sequence of $A$-modules, $M$ and $P$ be 
nonisomorphic simple modules and the base field $\K$ be algebraically closed.
Then the sequence $0\to\Fl M\to\Fl N$ is not exact.
\end{lemma}
\Proof We denote by $\rho_M\co A\to \End_{\K}M$ the representation of $A$ 
corresponding to the module $M.$ Choosing a $\K$-basis for $N,$ conforming 
to the composition series (\ref{eq1}), we obtain
$$
\begin{array}{c}
\rho_N(A)\subset
\left\{\left(
\begin{array}{cc}
\alpha &\beta\\
 0 &\gamma
\end{array}
\right)\right\} \text{ with $\alpha\in M_{n_1}(\K), 
\beta\in M_{n_1\times n_2}(\K), \gamma\in M_{n_2}(\K)$}\\
\text{ and } \rho_N(a)=\left(
\begin{array}{cc}
\alpha &\beta\\
 0 &\gamma
\end{array}
\right) \hence \rho_M(a)=\alpha,\ \rho_P(a)=\gamma.
\end{array}
$$
As $M\not\cong P$ are simple, we have $\rho_{M\oplus P}(A)=\End_{\K}M\oplus
\End_{\K}P,$ \ie $\rho_N(A)\myto{p}\rho_M(A)\oplus\rho_P(A)$ is an 
epimorphism; if it were an isomorphism, $\rho_N(A)$ would be semisimple 
and $N$ would be a direct sum of simple modules and a trivial extension, 
contradicting the hypothesis of the lemma. Thus
$$
\exists\beta\ne 0:
\left(\begin{array}{cc}
0 &\beta\\
 0 & 0
\end{array}
\right)
\in \rho_N(A);
$$ 
the matrix multiplication in $\rho_N(A)$ and the surjectivity of
$p$ show that $\beta$'s of this kind form a subrepresentation in the 
representation of $M_{n_1}(\K)\otimes M_{n_2}(\K)^0$ on
$M_{n_1\times n_2}(\K)$ by left and right multiplications respectively,
and the irreducibility of this representation shows that
$$
\rho_N(A)=
\left\{\left(
\begin{array}{cc}
\alpha &\beta\\
 0 &\gamma
\end{array}
\right)\right\} \text{ for $\alpha\in M_{n_1}(\K), 
\beta\in M_{n_1\times n_2}(\K), \gamma\in M_{n_2}(\K).$}
$$
So $$\Fl N=\T N\left/\left<\left(
\begin{array}{cc}
A &B\\
 0 & C
\end{array}
\right)\right>\right.,$$
$A,B,C$ being generic matrices of corresponding sizes, say,
$A=(a_{ij}),$ $B=(b_{ij}),$ $C=(c_{ij}),$
and $\Fl M=\T M/\langle A\rangle,$ $\Fl P=\T P/\langle
C\rangle.$ This shows that $\set{M}$, the image of $\Fl M$ in $\Fl N$,
equals
$$
\T M\left/\T M\cap\left<\left(
\begin{array}{cc}
A &B\\
 0 & C
\end{array}
\right)\right>\right..
$$
Now we show that
\begin{equation}\label{*}
\T M\cap\left<\left(
\begin{array}{cc}
A &B\\
 0 & C
\end{array}
\right)\right>\ne\langle A\rangle.
\end{equation}
Let $\sum M_j\eps_j\in\langle\eps_1,\dots,\eps_{ln}\rangle_R$ be a basis 
syzygy on the columns of the matrix $C.$ Then $\sum M_jb_j$ ($b_j$ 
being the $j$th column of $B$) belongs to the left-hand side of (\ref{*}) 
and depends only on $b_{\cdot\cdot}$ and $c_{\cdot\cdot},$ as in the basis 
syzygies on the columns of $C$ $M_j$ are $n\times n$-cofactors in some
$n\times(n+1)$ submatrix of $C$ (here we use that $l$ is greater than $1$ 
and there are such submatrices).

Thus the left-hand side of (\ref{*}) is strictly larger than the right-hand
side and $\set{M}\not=\Fl M,$ \ie the functor is not exact. \qed

So we see that every indecomposable $A$-module has only one type of 
composition factors. Since for a simple module $P$ the module $F_l(P)$ is 
always Cohen-Macaulay and the functor $F_l$ is exact, the behavior of the 
depth (Prop.~\ref{dpft}) and Krull dimension \cite[chap.~III, B.1, chap.~I,
C.1, prop.~10]{Se} in short exact sequences together with the induction 
on the length show that for every indecomposable $M$ $F_l(M)$ is 
Cohen-Macaulay. So the water is poured out of the kettle and we can apply
the remaining part of the theorem.

\chapter{Proof of Theorem~3}
{\bf Reduction to the module $P$.} First, we can pass to the algebraic
closure of the base field, as all these invariants are determined by the
structure of the minimal resolution of $\Fl M,$ so we assume that $\K$ is 
algebraically closed. Further, the Hilbert series
and the multiplicity are additive in short exact sequences in the category
of graded modules, so these invariants for $\Fl M$ are $\dim_{\K}M/\dim_{\K}P$
times as big as for $\Fl P,$ and $\dim_{\K}P=n.$ For Betti
numbers, in particular (Prop.~\ref{CMt}) for the Cohen-Macaulay type, this
follows from the fact that with Notations~\ref{Not} the minimal resolution
of $\Fl M$ equals $C\otimes_{R[T]}\T\M,$ so the Betti numbers are 
 $\rk_R\T\M=\dim_{\K}M/n$ times as big as the ranks of the components of
$C,$ and for $P$ this coefficient equals 1. Thus it suffices to calculate
all the invariants for $H_0(C).$ We do the calculations for an
Eagon-Northcott complex of a $(g\times f)$-matrix and then substitute our
values $f=ln, g=n.$

{\bf The Betti numbers of $\Fl P$.} 
The formulas are immediate from the definition of the Eagon-Northcott
complex: the matrix $\varphi$ has only entries of degree $1$, so, to give 
the differentials in the Eagon-Northcott complex degree zero when the
generators of $F_0$ have degree zero, the generators of $F_1$ should have 
degree $1,$ the generators of $F_2$ degree $g+1$, as the elements of the 
matrix of the differential are $g\times g$-minors of $\varphi,$ and then the 
degree should advance by $1.$ As the rank of $\bigwedge^i F\otimes(S^j 
G)^*$ equals $\binom{f}{i}\binom{g+j-1}{g-1}$, we get the required values 
for Betti numbers.

{\bf The Hilbert series and the multiplicity.} The Hilbert series is an
additive function on the graded modules, so the Euler-Poincare characteristic
of the minimal resolution
$$0\to F_k\to\dots\to F_0\to\Fl P\to 0$$ of the module $\Fl P$ \wrt the
Hilbert series equals zero, \ie $\Fl P(t)=\sum_i (-1)^i F_i(t).$ If we
consider the Poincare series
$P(s,t)=\sum_is^it^jb_{ij},$ $b_{ij}$ being the graded Betti numbers of 
$\Fl P,$ then $\Fl P(t)=P(-1,t)/(1-t)^{\dim R}.$ We have:

$$
P(s,t)=g+fst+\sum_{k=g+1}^{f}\binom{f}{k}\binom{k-2}{g-1} s^{k+1-g}t^k=
$$$$
=(1-(-s)^{1-g})(g+fst)+\frac1{(g-1)!}s^2
\left(\frac{\partial}{\partial s}\right)^{g-1}\left(s^{-2}(1+st)^f\right)
=
$$$$
=(1-(-s)^{1-g})(g+fst)+(1+st)^{f-g+1}
\sum_{k=0}^{g-1}\binom{f}{k}(g-k)t^k(-1/s-t)^{g-1-k}.
$$
(if we expand the $f$th power in the middle line by the binomial formula 
and then differentiate termwise, we obtain the previous expression, and if 
we differentiate the product $s^{-2}\cdot(1+st)^f$ according to the 
Leibniz formula, we obtain the last expression).

Substituting $s=-1$ and cancelling common factors with the denominator we
obtain the required expression:
$$
F_l(P)(t)=(1-t)^{f-g+1-\dim R}
\sum_{i=0}^{g-1}\binom{f}{i}(g-i)t^i(1-t)^{g-1-i}
$$
(here, we recall, $f=ln,$ $g=n$), and the sum is exactly the polynomial 
$p(t)$ in the definition of the multiplicity. Substituting $t=1$, we
find the multiplicity: only the summand with $i=g-1$ remains.\qed


\begin{thebibliography}{20}
\addcontentsline{toc}{chapter}{\bibname}
%
\bibitem{palam}
\textsc{Adams~W.~W., Loustaunau~P., Palamodov~V.~P., Strup\-pa~D.~C.}
Hartog's Phenomenon for Polyregular Functions and Projective Dimension of
Related Modules over a Polynomial Ring $//$ Ann. Inst. Fourier 1997, 47(2).
P.~623--640.
%
\bibitem{AL}
\textsc{Adams~W.~W., Loustaunau~P.} Analysis of the Module Determining the
Properties of Regular Functions of Several Quaternionic Variables $//$
Pacific J.~Math. 2000, 196(1). P.~1--15.
%
\bibitem{AzN}
\textsc{Azumaya\:G., Nakayama\:T.} On absolutely uni-serial algebras $//$ 
Jap. Journ. Math. 1948. V.~19. No.~4. P.~263--273.
%
\bibitem{Az}
\textsc{Azumaya\:G.} On maximally central algebras $//$ Nagoya Math. J. 1951.
V.~2. P.~119--150.
%
\bibitem{Bour}
\textsc{Bourbaki N.}
\'El\'ements de math\'ematiques. XXIII. Premi\'ere partie: Les structures
fondamentales de l'analyse. Livre II: Alg\`ebre. Chap. 8: Modules et anneaux
semisimples. Paris: Hermann, 1958. 
%
\bibitem{BourX}
\textsc{Bourbaki N.} \'El\'ements de math\'ematiques. Alg\`ebre. 
Chapitre 10: Alg\`ebre homologique. Paris: Masson, 1980. 
%
\bibitem{BH}
\textsc{Bruns\:W., Herzog\:J.} Cohen-Macaulay Rings. Cambridge Univ. 
Press, 1993. 
%
\bibitem{BV}
\textsc{Bruns\:W., Vetter\:U.} Determinantal Rings. Springer, 1988. (Lecture
notes in math. No. 1327.)
%
\bibitem{WMCE}
\textsc{Buchsbaum~D.~A., Eisenbud~D.} What Makes a Complex Exact? $//$ J.
Algebra 1973, 25(2). P.~259--268.
%
\bibitem{Maxw}
\textsc{Colombo\:F., Loustaunau\:P., Sabadini\:I., Struppa\:D.\:C.} Regular
functions of biquaternionic variables and Maxwell's equations // J. Geom.
Phys. 1998, 26(3--4). P.~183--201.
%
\bibitem{CR} 
\textsc{Curtis Ch., Reiner I.} Representation theory of finite groups and
associative algebras. New York-London: Interscience, 1962.
%
\bibitem{DeM}
\textsc{DeMeyer~F., Ingraham~E.} Separable Algebras over Commutative 
Rings. Springer, 1971. (Lecture notes in math. No. 181.)
%
\bibitem{Kro}
\textsc{Deuring\:M.} Algebren. 2.\@ Auflage. Springer, 1968. (Ergebnisse der 
Mathematik und ihrer Grenzgebiete. Bd. 41.) 
%
\bibitem{Ei}
\textsc{Eisenbud~D.} Commutative Algebra with a View Toward Algebraic
Geometry. Springer, 1995. (Graduate texts in math. No. 150.)
%
\bibitem{Her}
\textsc{Herstein\:I.\:N.}  Noncommutative rings. The Carus Mathematical
Monographs, No. 15. Published by The Mathematical Association of America;
distributed by John Wiley \& Sons, Inc., New York. xi+199 pp. (1968).
%
\bibitem{Ma}
\textsc{Manin Yu.~I.}
Cubic Forms. Algebra, Geometry, Arithmetic. North-Holland Mathematical Library.
Vol. 4. North-Holland, 1974.
%
\bibitem{Pal}
\textsc{Palamodov\:V.\:P.} Linear differential operators with constant
coefficients. Translated from the Russian by A. A. Brown. Die Grundlehren
der mathematischen Wissenschaften, Band 168. Berlin--Heidelberg--New York:
Springer-Verlag. VIII, 444 pp. (1970).
%
\bibitem{Moisil}
\textsc{Sabadini I., Shapiro M. V., Struppa D. C.}
Algebraic analysis of the Moisil-Theodorescu system // Complex Variables
Theory Appl. 2000, 40(4). P.~333--357.
%
\bibitem{Dirac}
\textsc{Sabadini I., Struppa D. C., Sommen~F., Van Lan\-cker~P.}
Complexes of Dirac operators in Clifford algebras //
Math. Z. 2002, 239(2). P.~293--320. 
%
\bibitem{Se}
\textsc{Serre~J.-P.} Alg\`ebre locale --- multiplicit\'es. Berlin:
Springer-Verlag, 1965. (Lecture notes in math. No. 11.)
%
\bibitem{VDW}
\textsc{van der Waerden~B.~L.}  
Algebra. I. Unter Benutzung von Vorlesungen von E. Artin und E. Noether. 8.
Auflage der Modernen Algebra. Heidelberger Taschenb\"ucher. Band 12.
Berlin--Heidelberg--New York: Springer-Verlag. IX, 272 S. (1971). 

Algebra. 2.\@ Teil. Unter Benutzung von Vorlesungen von E. Artin und E.
Noether. 5.\@ Aufl. der Modernen Algebra. Heidelberger Taschenb\"ucher. 23.
Berlin--Heidelberg--New York: Springer-Verlag. XII, 300 S. (1967).
%

\section*{Papers by the author}
%\addcontentsline{toc}{chapter}{\numberline{}Papers by the author}
%
\bibitem{art0}
\textsc{Popov\:O.\:N.} On a construction of modules over a polynomial 
ring. // Russian Mathematical Surveys, 56(6) (2001), 1177--1178.
%
\bibitem{art1} 
\textsc{Popov\:O.\:N.} Modules over a polynomial ring obtained
from representations of finite-dimensional associative algebras $//$ Sbornik:
Mathematics, 193(3) (2002), 423--443.
%
\bibitem{art2}
\textsc{Popov\:O.\:N.} More about a
construction for modules over a polynomial ring $//$ Russian Mathematical
Surveys, 58(2) (2003), 386--387.  
%
\bibitem{art3}
\textsc{Popov\:O.\:N.} 
On modules over a polynomial ring obtained from representations of
finite-dimensional associative algebras. II. The case of a non-perfect
field $//$ Sbornik: Mathematics, 195(9) (2004), 1309--1319. 
%
\bibitem{art4kurz}
\textsc{Popov\:O.\:N.} On a construction of modules over a
polynomial ring in the case of an arbitrary field $//$ Russian Mathematical
Surveys, 59(3) (2004), 583--584. 
%
\end{thebibliography}
\end{document}